%% file: maincasc.tex
\documentclass[envcountsame,envcountsect,runningheads,oribibl]{llncs}

\input{header}

\usepackage[pdfauthor   = {B\"achler and Gerdt and Lange-Hegermann and Robertz},
            pdftitle    = {Thomas Decomposition of Algebraic and Differential Systems},
            pdfsubject  = {},
            pdfkeywords = {},
            bookmarks=true,
            bookmarksopen=true,
            pagebackref=true,
            hyperindex=true,
            colorlinks=true,
            linkcolor=black,
            citecolor=black,
            filecolor=black,
            urlcolor=black,
            ps2pdf=true
            ]{hyperref}

\title{Thomas Decomposition of Algebraic and Differential Systems}

\newcommand{\articletype}{combined}
\newcommand{\splitarticle}{casc}
\input{splitrefs}

\begin{document}

\maketitle

\begin{abstract}
 In this paper we consider disjoint decomposition of algebraic and non-linear partial differential systems of equations and inequations into so-called simple subsystems.
 We exploit \textsc{Thomas} decomposition ideas and develop them into a new algorithm.
 For algebraic systems simplicity means triangularity, squarefreeness and non-vanishing initials.
 For differential systems the algorithm provides not only algebraic simplicity but also involutivity.
 The algorithm has been implemented in \textsc{Maple}.
\end{abstract}

\input{introduction}
\input{algebraic}
\input{differential}

\input{implementation}
\input{acknowledgements}

\newcommand{\etalchar}[1]{$^{#1}$}
\def\cprime{$'$} \def\cprime{$'$}
\providecommand{\bysame}{\leavevmode\hbox to3em{\hrulefill}\thinspace}
\providecommand{\MR}{\relax\ifhmode\unskip\space\fi MR }
\providecommand{\MRhref}[2]{%
  \href{http://www.ams.org/mathscinet-getitem?mr=#1}{#2}
}
\providecommand{\href}[2]{#2}

\end{document}

%% file: header.tex
\usepackage[utf8]{inputenc}
\usepackage{a4wide, amsmath, amssymb}
\usepackage{algorithmic}
\usepackage[active]{srcltx}
\usepackage{tikz}
\usepackage{float}
\floatstyle{ruled}
\usepackage{enumerate}

\newfloat{Decompose}{tbp}{lop}[section]

\usepackage{maple2e}
 \def\emptyline{\vspace{6pt}}

\usepackage[T1]{fontenc}

\author{Thomas B\"achler\inst{1} \email{thomas@momo.math.rwth-aachen.de},
 Vladimir Gerdt\inst{2} \email{gerdt@jinr.ru},
 Markus Lange-Hegermann\inst{1} \email{markus@momo.math.rwth-aachen.de}
 and Daniel Robertz\inst{1} \email{daniel@momo.math.rwth-aachen.de}
}
\authorrunning{Thomas B\"achler, Vladimir Gerdt, Markus Lange-Hegermann and Daniel Robertz}

\institute{
  Lehrstuhl B f\"ur Mathematik, RWTH-Aachen University, Germany
  \and Joint Institute for Nuclear Research, Dubna, Russia
}

\spnewtheorem{myalgorithm}[theorem]{Algorithm}{\bfseries}{\rmfamily}
\spnewtheorem{mydefinition}[theorem]{Definition}{\bfseries}{\rmfamily}
\spnewtheorem{myremark}[theorem]{Remark}{\bfseries}{\rmfamily}
\spnewtheorem{mycorollary}[theorem]{Corollary}{\bfseries}{\rmfamily}
\spnewtheorem{mylemma}[theorem]{Lemma}{\bfseries}{\rmfamily}
\spnewtheorem{defrem}[theorem]{Definition and Remark}{\bfseries}{\rmfamily}
\spnewtheorem{myexample}[theorem]{Example}{\bfseries}{\rmfamily}

\newcommand{\ld}{\operatorname{ld}}
\newcommand{\ini}{\operatorname{init}}
\newcommand{\rank}{\operatorname{rank}}
\newcommand{\prem}{\mathsf{prem}}
\newcommand{\pquo}{\mathsf{pquo}}
\newcommand{\sep}{\operatorname{sep}}

\newcommand{\res}{\operatorname{res}}

\newcommand{\lcm}{\ensuremath{\operatorname{lcm}}}

\newcommand{\C}{\ensuremath{\mathbb{C}}}

\newcommand{\R}{\ensuremath{\mathbb{R}}}
\newcommand{\Z}{\ensuremath{\mathbb{Z}}}

\newcommand{\sol}{\ensuremath{\mathfrak{Sol}}}

\newcommand{\prs}{\ensuremath{\operatorname{PRS}}}
\newcommand{\ma}{\ensuremath{\mathbf{a}}}
\newcommand{\pma}[1][]{\ensuremath{\phi_{\ifthenelse{\equal{#1}{}}{}{#1,}\mathbf{a}}}}
\def\osim{
  \mathrel{
    \mathchoice {
      \vcenter{
        \offinterlineskip
        \halign{
          \hfil $\displaystyle##$\hfil\cr\not\succ\cr\noalign{\vskip-3pt}\not\prec\cr
        }
      }
    }
    {
      \vcenter{
        \offinterlineskip
        \halign{
          \hfil$\textstyle##$\hfil\cr\not\succ\cr\noalign{\vskip-3pt}\not\prec\cr
        }
      }
    }
    {
      \vcenter{
        \offinterlineskip
        \halign{
          \hfil$\scriptstyle##$\hfil\cr\not\succ\cr\noalign{\vskip-2.4pt}\not\prec\cr
        }
      }
    }
    {
      \vcenter{
        \offinterlineskip
        \halign{
          \hfil$\scriptscriptstyle##$\hfil\cr\not\succ\cr\noalign{\vskip-0.9pt}\not\prec\cr
        }
      }
    }
  }
}

%% file: splitrefs.tex
\usepackage{boolexpr}
\usepackage{ifthen}

\newcommand{\customcrossref}[1]%
{%
\switch[\pdfstrcmp{#1}]%
\case{{section_differential}}\cite{BGLR2}%
\case{{exist_sol}}\cite[Rem.~2.2]{BGLR1}%
\case{{algebraic_definition_notation}}\cite[\textsection2.1]{BGLR1}%
\case{{algebraic_algorithms}}\cite[\textsection2.2]{BGLR1}%
\case{{algo_decompose}}\cite[Alg.~2.20]{BGLR1}%
\case{{prempquo}}\cite[Eq.~(1)]{BGLR1}%
\case{{reduce_q}}\cite[Alg.~2.20, line~8]{BGLR1}%
\case{{q_2}}\cite[Alg.~2.20, line~25]{BGLR1}%
\case{{resplitgcd}}\cite[Alg.~2.20, line~14]{BGLR1}%
\case{{insertres}}\cite[Alg.~2.20, line~17]{BGLR1}%
\case{{dividebyinequation}}\cite[Alg.~2.20, line~29]{BGLR1}%
\case{{q_3}}\cite[Alg.~2.20, line~30]{BGLR1}%
\case{{def_orders}}\cite[2.24]{BGLR1}%
\case{{composite_well_founded}}\cite[Rem.~2.23]{BGLR1}%
\case{{BGLROptimizations}}\cite[\textsection3.1]{BGLR1}%
\otherwise\textbf{??}%
\endswitch%
}

\newcommand{\crossref}[3]{%
 \ifthenelse{\equal{\articletype}{combined}}
 {#1\ref{#2}#3}
 {\customcrossref{#2}}%
}

%% file: introduction.tex
\section{Introduction}\label{introduction}

\ifthenelse{\equal{\articletype}{combined} \or \equal{\splitarticle}{algebraic}}{%
Nowadays, triangular decomposition algorithms, which go back to the characteristic set method by Ritt \cite{Ritt} and Wu \cite{Wu}, and software implementing them have become powerful tools for investigating and solving systems of multivariate polynomial equations.
In many cases these methods are computationally more efficient than those based on construction of \textsc{Gr\"obner} bases. As an example of such problems one can indicate \textsc{Bool}ean polynomial systems arising in cryptanalysis of stream ciphers.
For those systems triangular decomposition algorithms based on the characteristic set method revealed their superiority over the best modern algorithms for the construction of \textsc{Gr\"obner} bases \cite{GaoHuang}.

For terminology, literature, definitions and basic proofs on triangular-decomposition algorithms for polynomial and differential-polynomial systems we refer to the excellent tutorial papers
\cite{Hubert1,Hubert2} and to the bibliographical references therein.

Among numerous triangular decompositions the \textsc{Thomas} one stands by itself.
It was suggested by the American mathematician J.M.Thomas in his books \cite{Tho1,Tho2} and decomposes a finite system of polynomial equations and inequations into finitely many triangular subsystems that he called \emph{simple}.
Unlike other decomposition algorithms it yields a \emph{disjoint} zero decomposition, that is, it decomposes the affine variety or quasi-affine variety defined by the input into a finite number of disjoint quasi-affine varieties determined by the output simple systems.
Every simple system is a regular chain.
}{%
This article is a follow-up of \cite{BGLR1}, where the algebraic \textsc{Thomas} decomposition was treated in detail.
Here we extend the algebraic decomposition to differential systems.
The reader should be familiar with \cite{BGLR1}, since notations and technical details from there are used throughout this article extensively.
}%

\ifthenelse{\equal{\articletype}{combined} \or \equal{\splitarticle}{differential}}{%
During his research on triangular decomposition, Thomas was motivated by the \textsc{Riquier}-\textsc{Janet} theory \cite{Riquier,Janet}, extending it to non-linear systems of partial differential equations.
For this purpose he developed a theory of (\textsc{Thomas}) monomials, which generate the involutive monomial division called \textsc{Thomas} division in \cite{GB1}.
He gave a recipe for decomposing a non-linear differential system into algebraically simple and passive subsystems \cite{Tho1}.

Differential \textsc{Thomas} decomposition differs noticeably from that computed by the famous \textsc{Rosenfeld}-\textsc{Gr\"obner} algorithm \cite{BLOP,BLOP2} which forms a basis of the \textsf{diffalg} and \textsf{BLAD} libraries \cite{diffalg,blad} as well as from other differential decompositions (e.g. \cite{Bouziane}). 
We found that \textsf{diffalg} and \textsf{BLAD} are optimized and well-suited for ordinary differential equations.
However, all other known methods give a zero decomposition which, unlike that in \textsc{Thomas} decomposition, is not necessarily disjoint.
}{}%

\ifthenelse{\equal{\splitarticle}{jsc}}{%
Thus, similar to \textsf{diffalg}, the %
}{}%
\ifthenelse{\equal{\splitarticle}{algebraic}}{%
The %
}{}%
\ifthenelse{\equal{\splitarticle}{jsc} \or \equal{\splitarticle}{algebraic}}{%
\textsf{RegularChains} package \cite{regularchains}, which is shipped with recent versions of \textsc{Maple}, implements a decomposition of a polynomial ideal into ideals represented by regular chains and a radical decomposition of an ideal into square-free regular chains.
The solution sets of this decomposition are %
\ifthenelse{\equal{\articletype}{combined}}{%
also%
}{}%
not disjoint in general.
There is an extension called comprehensive triangular decomposition (cf.~\cite{ctd}) that provides disjointness on the parameters of a parametric system, but does in general not produce a decomposition into simple systems.
The disjointness of the \textsc{Thomas} decomposition combined with the properties of the simple systems provide a useful platform for counting solutions of polynomial systems.
In fact, the {\em counting polynomial} introduced by Plesken in \cite{PleskenCounting} can only be computed using a \textsc{Thomas} decomposition.
An extension of this polynomial to the differential case is in development.
}{}%

\ifthenelse{\equal{\splitarticle}{jsc} \or \equal{\splitarticle}{algebraic}}{%

It can be used to compute invariants of affine or projective varieties like their dimension.
}{}%
\ifthenelse{\equal{\splitarticle}{jsc}}{%
Further applications including counting of solutions in the differential case are under development.
}{}%

A first implementation of the \textsc{Thomas} decomposition was done by Teresa Gómez-Díaz in AXIOM under the name
``dynamic constructible closure'' which later turned out to be the same as the \textsc{Thomas}
decomposition \cite{delliere_wang_dyn}.
Wang later designed and implemented an algorithm constructing the \textsc{Thomas} decomposition \cite{wang_simple,WangMethods,WangLi}.
For polynomial and ordinary differential systems Wang's algorithm was implemented by himself in \textsc{Maple} \cite{WangPractice} as part of the software package $\epsilon$\textsf{psilon} \cite{epsilon}, which also contains implementations of a number of other triangular decomposition algorithms.
A modified algorithmic version of the \textsc{Thomas} decomposition was considered in \cite{GerdtSimple}%
\ifthenelse{\equal{\articletype}{combined} \or \equal{\splitarticle}{differential}}{%
\ with its link to the theory of involutive bases \cite{GB1,GerI,Ger3}%
}{}%
.
\ifthenelse{\equal{\articletype}{combined} \or \equal{\splitarticle}{differential}}{%
The latter theory together with some extensions is presented in detail in the recent book \cite{Seiler}.
}{}%

In the given paper we present a new algorithmic version of the \textsc{Thomas} decomposition for %
\ifthenelse{\equal{\articletype}{combined}}{%
polynomial and (partial) differential %
}{%
\ifthenelse{\equal{\splitarticle}{algebraic}}{%
polynomial %
}{%
differential %
}}%
systems.
\ifthenelse{\equal{\articletype}{combined}}{%
In the differential case %
}{}%
\ifthenelse{\equal{\splitarticle}{differential}}{%
In this case %
}{}%
\ifthenelse{\equal{\articletype}{combined} \or \equal{\splitarticle}{differential}}{%
the output subsystems are \textsc{Janet} involutive in accordance to the involutivity criterion from \cite{GerdtSimple}, and hence they are coherent.
Moreover, for every output subsystem the set of its equations is a minimal \textsc{Janet} basis of the radical differential ideal generated by this set.
}{}%
The algorithm has been implemented in \textsc{Maple}%
\ifthenelse{\equal{\articletype}{combined}}{%
\ for both the algebraic and differential case%
}{}%
.
\ifthenelse{\equal{\splitarticle}{algebraic}}{%
There is a follow-up \cite{BGLR2} to this article describing the differential case of the \textsc{Thomas} decomposition.
}{}%
\ifthenelse{\equal{\articletype}{combined} \or \equal{\splitarticle}{differential}}{%
For a linear differential system it constructs a \textsc{Janet} basis of the corresponding differential ideal and for this case works similarly to the \textsc{Maple} package \textsf{Janet} (cf.\ \cite{JanetPackage}).
}{}%

This paper is organized as follows.
\ifthenelse{\equal{\articletype}{combined} \or \equal{\splitarticle}{algebraic}}{%
In \textsection\ref{algebraic} we %
\ifthenelse{\equal{\splitarticle}{jsc}}{present}{sketch}
\ifthenelse{\equal{\articletype}{combined}}{%
the algebraic part of %
}{}%
our algorithm for the \textsc{Thomas} decomposition with its main objects defined in \textsection\ref{algebraic_definition_notation}.
The algorithm itself together with its subalgorithms is considered in \textsection\ref{algebraic_algorithms}\ifthenelse{\equal{\articletype}{split} \or \equal{\splitarticle}{jsc}}{, where we prove its correctness and termination}{}.
\ifthenelse{\equal{\splitarticle}{jsc}}{ %
The work of the algorithm is illustrated by an example from which one can trace the results of the algebraic transformations performed by the main steps of the algorithm. %
}
}{}%
\ifthenelse{\equal{\articletype}{combined}}{%
Decomposition of differential systems is described in \textsection\ref{section_differential}.
Here we %
}{}%
\ifthenelse{\equal{\splitarticle}{differential}}{%
We %
}{}%
\ifthenelse{\equal{\articletype}{combined} \or \equal{\splitarticle}{differential}}{%
briefly introduce some basic notions and concepts from differential algebra (\textsection\ref{differential_preliminaries}) and from the theory of involutive bases specific to \textsc{Janet} division (\textsection\ref{Janet}) together with one of the two extra subalgorithms that extend the algebraic decomposition %
\ifthenelse{\equal{\splitarticle}{differential}}{%
from \cite{BGLR1} %
}{}%
to the differential one.
The second such subalgorithm is considered in \textsection\ref{differential reduction} along with the definition of differential simple systems.
Subsection \textsection\ref{differential algorithm} contains a description of the differential \textsc{Thomas} decomposition algorithm\ifthenelse{\equal{\splitarticle}{jsc}}{and the proof of its correctness and termination}{}.
}{}%
Some implementation issues are discussed in \textsection\ref{implementation}%
\ifthenelse{\equal{\splitarticle}{algebraic}}{%
.%
}{%
, where we also demonstrate the \textsc{Maple} implementation for the differential decomposition using the example of a system related to control theory.
}%

\ifthenelse{\equal{\splitarticle}{casc}}{We omit the proofs for compactness. They will be published elsewhere.}{}

%% file: algebraic.tex
\section{Algebraic Systems}\label{algebraic}

The algebraic \textsc{Thomas} decomposition deals with systems of polynomial equations and inequations.
This section introduces the concepts of simple systems and disjoint decompositions based on properties of the set of solutions of a system.
A pseudo reduction procedure and several splitting algorithms on the basis of polynomial remainder sequences are introduced as tools for the main algorithm, which is presented at the end of the section.

\subsection{Preliminaries}\label{algebraic_definition_notation}

Let $F$ be a computable field of characteristic 0 and $R:=F[x_1,\dots,x_n]$ the polynomial ring in $n$ variables.
A total order $<$ on the indeterminates of $R$ is called a \textbf{ranking}. The notation $R=F[x_1,\dots,x_n]$ shall implicitly define the ranking $x_1 < \ldots < x_n$.
The indeterminate $x$ is called \textbf{leader} of $p \in R$ if $x$ is the $<$-largest variable occurring in $p$ and we write $\ld(p)=x$.
If $p \in F$, we define $\ld(p)=1$ and $1<x$ for all indeterminates $x$.
The degree of $p$ in $\ld(p)$ is called \textbf{rank} of $p$ and the leading coefficient $\ini(p) \in F[\ y\ |\ y<\ld(p)\ ]$ of $\ld(p)^{\rank(p)}$ in $p$ is called \textbf{initial} of $p$.

For $\ma \in \overline{F}^n$, where $\overline{F}$ denotes the algebraic closure of $F$, define the following evaluation homomorphisms:
\[\pma: F[x_1,\dots,x_n] \to \overline{F}: x_i \mapsto a_i\]
\[\pma[<x_k]: F[x_1,\ldots,x_n] \to \overline{F}[x_k,\ldots,x_n]: \left\{\begin{array}{ll} x_i \mapsto a_i,&i<k \\x_i \mapsto x_i,&\mbox{otherwise}\end{array}\right.\]

For a polynomial $p \in R$, the symbols $p_{=}$ and $p_{\neq}$ shall denote the equation $p=0$ and inequation $p\neq0$, respectively.
A finite set of equations and inequations is called an \textbf{(algebraic) system} over $R$.
Abusing notation, we sometimes treat $p_{=}$ or $p_{\neq}$ as the underlying polynomial $p$.
A \textbf{solution} of a system $S$ is a tuple $\ma \in \overline{F}^n$ such that $\pma(p)=0$ for all equations $p_{=} \in S$ and $\pma(p)\neq0$ for all inequations $p_{\neq} \in S$.
The set of all solutions of $S$ is denoted by $\sol(S)$.

Define $S_x := \{ p \in S\ |\ \ld(p)=x \}$.
In a situation where it is clear that $|S_x|=1$, we also use $S_x$ to denote the unique element of $S_x$.
The subset $S_{<x} := \{ p \in S\ |\ \ld(p)<x \}$ can be considered a system over $F[\ y\ |\ y<x\ ]$.
Furthermore, the sets of all equations $p_{=} \in S$ and all inequations $p_{\neq}\in S$ are denoted by $S^{=}$ and $S^{\neq}$, respectively.

The general idea of the {\sc Thomas} methods is to use the homomorphism $\pma[<x]$ to treat each polynomial $p\in S_x$ as the \emph{univariate} polynomial $\pma[<x](p) \in \overline{F}[x]$ for all $\ma \in \sol(S_{<x})$ \emph{simultaneously}.
This idea forms the basis of our central object, the \textbf{simple system}:

\begin{mydefinition}[Simple Systems]\label{simple_system} Let $S$ be a system.
 \begin{enumerate}
  \item $S$ is \textbf{triangular} if $|S_{x_i}|\leq 1\ \forall\ 1\leq i \leq n$ and $S \cap \{ c_=, c_{\neq}\mid c \in F\}=\emptyset$.
  \item \label{simple_nonzeroinitials} $S$ has \textbf{non-vanishing initials} if $\pma(\ini(p)) \not=0\ \forall\ \ma \in \sol(S_{<x_i})$ and $p \in S_{x_i}$ for $1 \leq i \leq n$.
  \item $S$ is \textbf{square-free} if the univariate polynomial $\pma[<x_i](p) \in \overline{F}[x_i]$ is square-free $\forall\ \ma \in \sol(S_{<x_i})$ and $p \in S_{x_i}$ for  $1\leq i\leq n$.
  \item $S$ is called \textbf{simple} if it is \emph{triangular}, has \emph{non-vanishing initials} and is \emph{square-free}.
 \end{enumerate}
\end{mydefinition}

Although all required properties are characterized via solutions of lower-ranking equations and inequations, the {\sc Thomas} decomposition algorithm does not calculate solutions of polynomials.
Instead, it uses polynomial equations and inequations to \emph{partition} the set of solutions of the lower-ranking system to ensure the above properties.
\begin{myremark}\label{exist_sol}
Simplicity of a system guarantees the existence of solutions:
If $\mathbf{b} \in \sol(S_{<x})$ and $S_x$ is not empty, then $\phi_{<x,\mathbf{b}}(S_x)$ is a univariate polynomial with exactly $\rank(S_x)$ \emph{distinct} roots.
When extending $\mathbf{b}$ to a solution $(\mathbf{b},a)$ of $S_{\le x}$,
for an equation in $S_x$ there are $\rank(S_x)$ choices for $a$, whereas for an inequation or empty $S_x$ all but finitely many $a \in \overline{F}$ give an extension.
\end{myremark}

To transform a system into a simple system, it is in general necessary to partition the set of solutions.
Instead of an equivalent simple system, this leads to a so-called decomposition into simple systems.

\begin{mydefinition}
 A family $(S_i)_{i=1}^m$ is called \textbf{decomposition} of $S$ if $\sol(S)=\bigcup_{i=1}^m \sol(S_i)$.
 A decomposition is called \textbf{disjoint} if $\sol(S_i) \cap \sol(S_j) = \emptyset\ \forall\ i\neq j$.
 A \emph{disjoint} decomposition of a system into \emph{simple systems} is called \textbf{(algebraic) \textsc{Thomas} decomposition}.
\end{mydefinition}

For any algebraic system $S$, there exists a {\sc Thomas} decomposition (cf.\ \cite{Tho1}, \cite{Tho2}, \cite{wang_simple}).
The algorithm presented in the following section provides another proof of this fact.
First, we give an easy example of a \textsc{Thomas} decomposition.

\begin{myexample}Consider the equation\\
\begin{tabular}{cc}  
  \begin{minipage}{0.6\textwidth}
  \[
    p=y^2-x^3-x^2\enspace.
  \]
  A {\sc Thomas} decomposition of $\{p_=\}$ is given by:
  \[
    \left(\{(y^2-x^3-x^2)_=, (x\cdot(x+1))_{\not=} \}, \{y_=, (x\cdot(x+1))_= \}\right)
  \]
  \end{minipage}
  
  &
  
  \begin{minipage}{0.3\textwidth}
    \input{curve}
  \end{minipage}
\end{tabular}
\end{myexample}

\subsection{Decomposition Algorithms}\label{algebraic_algorithms}

Our version of the decomposition algorithm in each round treats one system, potentially splitting it into several subsystems. For this purpose, one polynomial is chosen from a list of polynomials to be processed.
This polynomial is pseudo-reduced modulo the system and afterwards combined with the polynomial in the system having the same leader.
To ensure that all polynomials are square-free and their initials do not vanish, the system may be split into several ones by initials of polynomials or subresultants.

From now on, a system $S$ is presented as a pair of sets $(S_T, S_Q)$, where $S_T$ represents a candidate for a simple system while $S_Q$ is the queue of elements to be processed.
$S_T$ is always triangular and $(S_T)_x$ denotes the unique equation or inequation of leader $x$ in $S_T$, if any.
$S_T$ also fulfills a weaker form of the other two simplicity conditions, i.e., for any solution $\ma$ of $(S_T)_{<x} \cup (S_Q)_{<x}$, we have $\pma(\ini((S_T)_x)) \neq 0$ and $\pma[<x]((S_T)_x)$ is square-free.

From now on, let $\prem$ be a \textbf{pseudo remainder algorithm}\footnote{In our context $\prem$ does not necessarily have to be the classical pseudo remainder, but any sparse pseudo remainder with property (\ref{prempquo}) will suffice.} in $R$ and $\pquo$ the corresponding \textbf{pseudo quotient algorithm}, i.e., for $p$ and $q$ with $\ld(p)=\ld(q)=x$ \begin{equation}\label{prempquo}m\cdot p = \pquo(p,q,x) \cdot q+\prem(p,q,x)\end{equation} where $\deg_x(q)>\deg_x(\prem(p,q,x))$ and $m \in R\setminus\{0\}$ with $\ld(m)<x$ and $m \mid \ini(q)^k$ for some $k \in \mathbb{Z}_{\geq0}$.
Note that if the initials of $p$ and $q$ are non-zero, the initial of $\pquo(p,q,x)$ is also non-zero.
Equation (\ref{prempquo}) only allows us to replace $p$ with $\prem(p,q,x)$ if $m$ does not vanish on any solution.
The below Algorithm (\ref{algo_reduce}) and Remark (\ref{reduce0}) require the last property, which, by definition, holds in simple systems.

The following algorithm employs pseudo remainders and the triangular structure to reduce a polynomial modulo $S_T$:

\begin{myalgorithm}[{\sf Reduce}]\label{algo_reduce}\ \\
 \textit{Input:} A system $S$, a polynomial $p \in R$ \\
 \textit{Output:} A polynomial $q$ with $\pma(p)=0$ if and only if $\pma(q)=0$ for each $\ma\in\sol(S)$.\\
 \textit{Algorithm:} \begin{algorithmic}[1]
                      \STATE $x \gets \ld(p)$; $q\gets p$ 
                      \WHILE{$x>1$ and $(S_T)_x$ is an equation and $\rank(q) \geq \rank((S_T)_x)$}
                       \STATE $q \gets \prem(q, (S_T)_x, x)$
                       \STATE $x \gets \ld(q)$
                      \ENDWHILE
                      \IF{$x>1$ and $\textrm{\sf Reduce}(S, \ini(q)) = 0$}
                       \RETURN $\textrm{\sf Reduce}(S, q-\ini(q)x^{\rank(q)})$
                      \ELSE
                       \RETURN $q$
                      \ENDIF
                     \end{algorithmic}
\end{myalgorithm}
A polynomial $p$ is called \textbf{reduced modulo $S_T$} if $\mathsf{Reduce}(S, p)=p$.
A polynomial $p$ \textbf{reduces to $q$ modulo $S_T$} if $\mathsf{Reduce}(S, p)=q$.
\ifthenelse{\equal{\articletype}{split} \or \equal{\splitarticle}{jsc}}{
 \begin{proof}[Correctness]
  There exist $m \in R\setminus\{0\}$ with $\ld(m)<\ld(p)$ and $\pma(m) \neq 0$ for all $\ma \in \sol(S_{\le\ld(p)})$ such that
  \[
    \mathsf{Reduce}(S, p)=mp-\sum_{y \le x} c_y\cdot (S_T)_y
  \]
  with $c_y \in R$ and $\ld(c_y)\le \ld(p)$. This implies
  \[
    \pma(\mathsf{Reduce}(S, p)) = \underbrace{\pma(m)}_{\neq0}\pma(p)-\sum_{y \le x} \pma(c_y)\underbrace{\pma((S_T)_y)}_{=0}
  \]
  and therefore $\pma(p) = 0$ if and only if $\pma(\mathsf{Reduce}(S, p)) = 0$.\qed
 \end{proof}
}{}

The result of the {\sf Reduce} algorithm does not need to be a canonical normal form. It only needs to detect polynomials that vanish on all solutions of a system:

\begin{myremark}\label{reduce0}
 Let $p \in R$ with $\ld(p)=x$. $\mathsf{Reduce}(S, p)=0$ implies $\pma(p)=0\ \forall\ \ma\in\sol(S_{\leq x})$.
\end{myremark}
\ifthenelse{\equal{\articletype}{split} \or \equal{\splitarticle}{jsc}}{
 \begin{proof}
  For all $\ma \in \sol(S_{\le x})$, it holds that $\pma(p) = 0$ if and only if $\pma(\mathsf{Reduce}(S, p)) = 0$. The statement follows from $\pma(\mathsf{Reduce}(S, p)) = \pma(0) = 0$.\qed
 \end{proof}
}{}

The converse of this remark only holds if $(S_Q)_{\le x}=\emptyset$, i.e., $(S_T)_{\le x}$ is simple.
If it is not simple, but $\ld(p)=x$ and $(S_Q)^=_{<x}=\emptyset$ hold, we still have some information.
In particular, $\textsf{Reduce}(S,p)\neq0$ implies that either $\sol(S_{<x})$ is empty or there exists $\ma\in\sol(S_{<x} \cup \{ (S_T)_x \})$ such that $\pma(p)\neq0$.

We now direct our attention to the methods we use to produce disjoint decompositions.
Since $\left( S \cup \left\{ p_{\neq} \right\}, S \cup \left\{ p_{=} \right\} \right)$ is a disjoint decomposition of $S$, we will use the following one-line subalgorithm as the basis of all the splitting algorithms described below.

\begin{myalgorithm}[{\sf Split}]
 \textit{Input:} A system $S$, a polynomial $p \in R$ \\
 \textit{Output:} The disjoint decomposition $\left( S \cup \left\{ p_{\neq} \right\}, S \cup \left\{ p_{=} \right\} \right)$ of $S$.\\
 \textit{Algorithm:} \begin{algorithmic}[1]
             \RETURN $\left( \left(S_T, S_Q \cup \{ p_{\neq} \}\right), \left(S_T, S_Q \cup \{ p_{=} \}\right) \right)$
            \end{algorithmic}
\end{myalgorithm}

The output of the following splitting algorithms is not yet a disjoint decomposition of the input.
However, the main algorithm \textsf{Decompose} will use this output to construct a disjoint decomposition.
We single out these algorithms to make the main algorithm more compact and readable.
For details we refer to the input and output specifications of the algorithms in question.

The algorithm \textsf{InitSplit} ensures that in one of the returned systems the property \ref{simple_nonzeroinitials} in Definition (\ref{simple_system}) holds for the input polynomial.
In the other system the initial of that polynomial vanishes.

\begin{myalgorithm}[\sf InitSplit]\label{algo_initsplit}
 \textit{Input:} A system $S$, an equation or inequation $q$ with $\ld(q)=x$. \\
 \textit{Output:} Two systems $S_1$ and $S_2$, where $\left(S_1 \cup \{q\}, S_2\right)$ is a disjoint decomposition of $S \cup \{q\}$. Moreover, $\pma(\ini(q))\neq0$ holds for all $\ma \in \sol(S_1)$ and $\pma(\ini(q))=0$ for all $\ma \in \sol(S_2)$.
 \textit{Algorithm:} \begin{algorithmic}[1]
    \STATE $(S_1, S_2) \gets \textrm{\sf Split}(S, \ini(q))$
    \IF{$q$ is an equation}
     \STATE $(S_2)_Q \gets (S_2)_Q \cup \left\{ \left(q-\ini(q)x^{\rank(q)}\right)_{=} \right\}$
    \ELSIF{$q$ is an inequation}
     \STATE $(S_2)_Q \gets (S_2)_Q \cup \left\{ \left(q-\ini(q)x^{\rank(q)}\right)_{\neq} \right\}$
    \ENDIF
    \RETURN $(S_1, S_2)$
   \end{algorithmic}
\end{myalgorithm}

In Definition (\ref{simple_system}) we view a multivariate polynomial $p$ as the univariate polynomial $\pma[<\ld(p)](p)$.
For ensuring triangularity and square-freeness, we often compute the gcd of two polynomials, which generally depends on the inserted value $\mathbf{a}$.
Subresultants provide a generalization of the \textsc{Euclid}ean algorithm useful in our context and their initials distinguish the cases of different degrees of gcds.

\begin{mydefinition}\label{not_prs}
 Let $p, q \in R$ with $\ld(p)=\ld(q)=x$, $\deg_x(p)=d_p > \deg_x(q)=d_q$.
 We denote by $\prs(p,q,x)$ the \textbf{subresultant polynomial remainder sequence} (see \cite{habicht}, \cite[Chap.~7]{mishra}, \cite[Chap.~3]{Yap}) of $p$ and $q$ w.r.t. $x$, and  by $\prs_i(p,q,x)$, $i<d_q$ the regular polynomial of degree $i$ in $\prs(p,q,x)$ if it exists, or $0$ otherwise.
 Furthermore, $\prs_{d_p}(p,q,x):=p$, $\prs_{d_q}(p,q,x):=q$ and $\prs_i(p,q,x):=0$, $d_q<i<d_p$.

 Define $\res_i(p,q,x) := \ini\left(\prs_i\left(p,q,x\right)\right)$ for $0<i<d_p$, whereas $\res_{d_p}(p,q,x):=1$ and $\res_0(p,q,x):=\prs_0\left(p,q,x\right)$.
 Note that $\res(p,q,x):=\res_0(p,q,x)$ is the usual resultant.
\end{mydefinition}

\ifthenelse{\equal{\articletype}{split} \or \equal{\splitarticle}{jsc}}{
Our definitions are slightly different from the ones cited in the literature (\cite[Chap.~7]{mishra}, \cite[Chap.~3]{Yap}), since we only use the regular polynomials.
However, it is easy to see that all theorems from \cite[Chap.~7]{mishra} we refer to still hold for $i<d_q$.
}

\begin{mydefinition}\label{def_fibrcard}
 Let $S$ be a system and $p_1, p_2 \in R$ with $\ld(p_1)=\ld(p_2)=x$. If $|\sol(S_{<x})|>0$, we call
 \[
   i := \min\left\{ i \in \mathbb{Z}_{\ge0} \mid \exists\ \ma \in \sol(S_{<x}) \mbox{ such that } \deg_x(\gcd(\pma[<x](p_1),\pma[<x](p_2))) = i \right\}
 \]
 the \textbf{fiber cardinality} of $p_1$ and $p_2$ w.r.t. $S$. Moreover, if $(S_Q)_{<x}^==\emptyset$, then
 \[
   i^\prime := \min \{ i \in \mathbb{Z}_{\ge0} \mid \textrm{\sf Reduce}(\res_j(p_1,p_2,x), S_T)=0\ \forall\ j<i \textrm{ and } \textrm{\sf Reduce}(\res_i(p_1,p_2,x), S_T)\neq0\}
 \]
 is the \textbf{quasi fiber cardinality} of $p_1$ and $p_2$ w.r.t. $S$. A disjoint decomposition $(S_1, S_2)$ of $S$ such that \begin{enumerate}
   \item $\deg_x(\gcd(\pma[<x](p_1),\pma[<x](p_2))) = i\ \forall\ \ma \in \sol\left((S_1)_{<x}\right)$
   \item $\deg_x(\gcd(\pma[<x](p_1),\pma[<x](p_2))) > i\ \forall\ \ma \in \sol\left((S_2)_{<x}\right)$
  \end{enumerate}
 is called the $i$-th \textbf{fibration split} of $p_1$ and $p_2$ w.r.t. $S$. A polynomial $r \in R$ with $\ld(r)=x$ such that $\deg_x(r)=i$ and
 \[\pma[<x](r) \sim \gcd(\pma[<x](p_1),\pma[<x](p_2))\ \forall\ \ma \in \sol\left((S_1)_{<x}\right)\] is called the
 $i$-th \textbf{conditional greatest common divisor} of $p_1$ and $p_2$ w.r.t. $S$, where $p \sim q$ if and only if $p \in \overline{K}^*q$. Furthermore, $q \in R$ with $\ld(q)=x$ and $\deg_x(q)=\deg_x(p_1)-i$ such that
 \[\pma[<x](q) \sim \frac{\pma[<x](p_1)}{\gcd(\pma[<x](p_1),\pma[<x](p_2))}\ \forall\ \ma \in \sol\left((S_1)_{<x}\right)\] is called the $i$-th \textbf{conditional quotient} of $p_1$ by $p_2$ w.r.t. $S$.
 By replacing $\pma[<x](p_2)$ in the above definition with $\frac{\partial}{\partial x}(\pma[<x](p_1))$, we get an $i$-th \textbf{square-free split} and $i$-th \textbf{conditional square-free part} of $p_1$ w.r.t. $S$.
\end{mydefinition}

The fiber cardinality is often not immediately available, as we may be unable to take inequations into account.
However, we can use all information contained in the equations using reduction, if all equations are contained in $S_T$.
Thus we require $(S_Q)_{<x}^==\emptyset$ before doing any reduction.

In this situation, the quasi fiber cardinality is easy to calculate and in many cases will be identical to the fiber cardinality.
Furthermore, if we consider the system $S_2$ from an $i$-th fibration split of some polynomials for a system $S$ and ensure that $((S_2)_Q)_{<x}^==\emptyset$, then the quasi fiber cardinality of the same polynomials for $S_2$ will be $i+1$.
Therefore and due to the following lemma, the quasi fiber cardinality is good enough for our purposes.

\begin{mylemma}\label{lm_ressplit}
 Let $|\sol(S_{<x})|>0$ and $(S_Q)_{<x}^==\emptyset$. For $p_1$, $p_2$ as in Definition (\ref{def_fibrcard}) with $\pma(\ini(p_1))\neq0\ \forall\ \ma\in \sol(S_{<x})$ and $\rank(p_1)>\rank(p_2)$, let $i$ be the fiber cardinality of $p_1$ and $p_2$ w.r.t. $S$ and $i^\prime$ the corresponding quasi fiber cardinality.
 Then \[i^\prime \leq i\] where the equality holds if and only if $\left|\sol\left(S_{<x} \cup \{\res_{i^\prime}(p_1,p_2,x)_{\neq}\}\right)\right|>0$.
\end{mylemma}
\ifthenelse{\equal{\articletype}{split} \or \equal{\splitarticle}{jsc}}{
 \begin{proof}
 Let $\ma \in \sol(S_{<x})$, $\rank(p_1)>\rank(p_2)$, $d_{p_1} := \deg_x(p_1) = \deg_x(\pma[<x](p_1))$, $d_{p_2} := \deg_x(p_2)$ and $d_{p_2,\ma} := \deg_x(\pma[<x](p_2))$: If $i < \max(d_{p_1},d_{p_2,\ma})-1 = d_{p_1}-1$, then \cite[Thm.~7.8.1]{mishra} implies the conditions \begin{equation}\label{lm_ressplit_proof1}\pma[<x](\prs_i(p_1,p_2,x))=0 \Longleftrightarrow \prs_i(\pma[<x](p_1),\pma[<x](p_2),x)=0\end{equation} and \begin{equation}\label{lm_ressplit_proof2}\pma(\res_i(p_2,p_2,x))=0 \Longleftrightarrow \res_i(\pma[<x](p_1),\pma[<x](p_2),x)=0\enspace.\end{equation}

 As $\prs_{d_{p_2}}(p_1,p_2,x)=p_2$, $\prs_{d_{p_1}}(p_1,p_2,x)=p_1$ and $\prs_{j}(p_1,p_2,x)=0$, $d_{p_2}<j<d_{p_1}$, conditions (\ref{lm_ressplit_proof1}) and (\ref{lm_ressplit_proof2}) also hold for $d_{p_2} \le i \le d_{p_1}$.

 Due to $\mathsf{Reduce}(\res_j(p_1,p_2,x), S_T)=0$ for all $j<i^\prime$, Remark (\ref{reduce0}) implies $\pma(\res_j(p_1,p_2,x))=0$, which by (\ref{lm_ressplit_proof1}) and (\ref{lm_ressplit_proof2}) implies $\res_j(\pma[<x](p_1),\pma[<x](p_2),x)=0$.
 By applying \cite[Thm.~7.10.5]{mishra} successively, we get $\prs_j(\pma[<x](p_1),\pma[<x](p_2),x)=0$ and thus \begin{equation} \label{lm_ressplit_proof3} \deg_x(\gcd(\pma[<x](p_1),\pma[<x](p_2)))\ge i^\prime\end{equation} holds because those $\prs_{j^\prime}(\pma[<x](p_1),\pma[<x](p_2),x), 0\le j^\prime < d_{p_2,\ma}$ which are non-zero are the polynomials constructed by the \textsc{Euclid}ean algorithm for the univariate polynomials $\pma[<x](p_1)$ and $\pma[<x](p_2)$ up to similarity.
 This implies $i^\prime \le i$.

 There exists $\mathbf{b} \in \sol(S_{<x})$ such that $\phi_\mathbf{b}(\res_{i^\prime}(p_1,p_2,x))\neq0$ and thus equality holds in (\ref{lm_ressplit_proof3}) for $\ma=\mathbf{b}$ if and only if $\left|\sol\left(S_{<x} \cup \{\res_{i^\prime}(p_1,p_2,x)_{\neq}\}\right)\right|>0$.\qed
 \end{proof}
}{}

\ifthenelse{\equal{\articletype}{split} \or \equal{\splitarticle}{jsc}}{
As the above lemma only applies if the degree of the first polynomial is strictly larger than the degree of the second one, we need another statement to cover the remaining cases.
Before, only the first polynomial needed a non-vanishing initial, whereas in the following version, this must also hold for the second polynomial.
}

\begin{mycorollary}\label{cor_ressplit}
 Let $|\sol(S_{<x})|>0$ and $(S_Q)_{<x}^==\emptyset$. For polynomials $p_1$, $p_2$ as in Definition (\ref{def_fibrcard}) with $\pma(\ini(p_1))\neq0$ and $\pma(\ini(p_2))\neq0\ \forall\ \ma\in \sol(S_{<x})$, let $i$ be the fiber cardinality of $p_1$ and $p_2$ w.r.t. $S$ and $i^\prime$ the quasi fiber cardinality of $p_1$ and $\prem(p_2,p_1,x)$ w.r.t. $S$. Then \[i^\prime \leq i\] with equality if and only if $\left|\sol\left(S_{<x} \cup \{\res_{i^\prime}(p_1,\prem(p_2,p_1,x),x)_{\neq}\}\right)\right|>0$.
\end{mycorollary}
\ifthenelse{\equal{\articletype}{split} \or \equal{\splitarticle}{jsc}}{
 \begin{proof}
  Let $\ma\in \sol(S_{<x})$. Due to $\pma(\ini(p_1))\neq0\neq\pma(\ini(p_2))$, \cite[Corr.~7.5.6]{mishra} implies that $\pma[<x](\prem(p_2,p_1,x))=\prem(\pma[<x](p_2),\pma[<x](p_1),x)$. The pairs of univariate polynomials $(\pma[<x](p_1), \pma[<x](p_2))$ and $(\pma[<x](p_1), \prem(\pma[<x](p_2),\pma[<x](p_1),x))$ have the same greatest common divisor, and thus we can replace $p_2$ with $\prem(p_2,p_1,x)$ in Lemma (\ref{lm_ressplit}).\qed
 \end{proof}
}{}

The following algorithm calculates the quasi fiber cardinality of two polynomials.
It is used as the basis for all algorithms that calculate a greatest common divisor or a least common multiple.

\begin{myalgorithm}[\sf ResSplit]\label{algo_ressplit}
 \textit{Input:} A system $S$ with $(S_Q)_{<x}^==\emptyset$, two polynomials $p, q \in R$ with $\ld(p)=\ld(q)=x$, $\rank(p)>\rank(q)$ and $\pma(\ini(p))\neq0$ for all $\ma \in \sol(S_{<x})$.\\
 \textit{Output:} The quasi fiber cardinality $i$ of $p$ and $q$ w.r.t. $S$ and an $i$-th fibration split $(S_1,S_2)$ of $p$ and $q$ w.r.t. $S$.\\
 \textit{Algorithm:} \begin{algorithmic}[1]
    \STATE \label{ressplit_cond_nontrivial_split} $i \gets \min \{ i \in \mathbb{Z}_{\ge0} \mid \textrm{\sf Reduce}(\res_j(p,q,x), S_T)=0\ \forall\ j<i \textrm{ and } \textrm{\sf Reduce}(\res_i(p,q,x), S_T)\neq0$\}
    \RETURN $(i, S_1, S_2) := \left(i, \textrm{\sf Split}(S, \res_i(p,q,x)) \right)$
   \end{algorithmic}
\end{myalgorithm}
\ifthenelse{\equal{\articletype}{split} \or \equal{\splitarticle}{jsc}}{
\begin{proof}[Correctness]
 If $|\sol((S_i)_{<x})|=0$ for $i\in \{1,2\}$, the statement is trivial. Thus, let $|\sol((S_i)_{<x})|>0$, $i=1,2$. Due to $\sol((S_i)_{<x}) \subseteq \sol(S_{<x})$, the statements from the proofs of Lemma (\ref{lm_ressplit}) and Corollary (\ref{cor_ressplit}) still hold for $\ma \in \sol((S_i)_{<x})$, $i \in \{1,2\}$.

 Let $\ma \in \sol((S_1)_{<x})$.
 The polynomial $g:=\prs_i(\pma[<x](p),\pma[<x](q),x)\not\equiv0$, due to $(\ini(g))_{\not=}=(\res_i(p,q,x))_{\neq} \in (S_1)_Q$.
 The degree of $g$ is $i$ and $g$ is similar to the gcd of $\pma[<x](p)$ and $\pma[<x](q)$ as discussed in the proof of Lemma (\ref{lm_ressplit}).

 For $\ma \in\sol((S_2)_{<x})$, we conclude from $(\ini(g))_==(\res_i(p,q,x))_= \in (S_2)_Q$ and \cite[Thm.~7.10.5]{mishra} that $g\equiv0$ and therefore $\deg_x(\gcd(\pma[<x](p),\pma[<x](q))) > i$.\qed
\end{proof}
}{}

Similarly to the \textsf{InitSplit} algorithm (\ref{algo_initsplit}), the following algorithm does not return a disjoint decomposition, but \textsf{Decompose} uses it to construct one.
\ifthenelse{\equal{\articletype}{split} \or \equal{\splitarticle}{jsc}}{
Like the other following algorithms, it is based on \textsf{ResSplit}.
It is used to calculate a greatest common divisor of an existing polynomial in $S_T$ and a new polynomial.
As seen in Definitions (\ref{not_prs}) and (\ref{def_fibrcard}), we produce one system where we can use the returned polynomial as a gcd, and another one where the degree of the gcd will be higher.
The algorithm will be applied to the latter system again later.
}

\begin{myalgorithm}[\sf ResSplitGCD]\label{algo_ressplitgcd}
 \textit{Input:} A system $S$ with $(S_Q)_{<x}^==\emptyset$, where $(S_T)_x$ is an equation, and an equation $q_=$ with $\ld(q)=x$. Furthermore $\rank(q)<\rank((S_T)_x)$.\\
 \textit{Output:} Two systems $S_1$ and $S_2$ and an equation $\widetilde{q}_=$ such that:
   \renewcommand{\theenumi}{\alph{enumi}}
   \renewcommand{\labelenumi}{\theenumi)}
   \begin{enumerate}
    \item\label{algo_ressplitgcd_cond1} $S_2 = \widetilde{S_2} \cup \{ q \}$ where $\left(S_1, \widetilde{S_2}\right)$ is an $i$-th fibration split of $(S_T)_x$ and $q$ w.r.t. $S$
    \item\label{algo_ressplitgcd_cond2} $\widetilde{q}$ is an $i$-th conditional gcd of $(S_T)_x$ and $q$ w.r.t. $S$.
   \end{enumerate}
   \renewcommand{\theenumi}{\arabic{enumi}}
   \renewcommand{\labelenumi}{\theenumi)}
   where $i$ is the quasi fiber cardinality of $p$ and $q$ w.r.t. $S$.\\
 \textit{Algorithm:} \begin{algorithmic}[1]
      \STATE $(i, S_1, S_2) \gets \textrm{\sf ResSplit}\left(S, (S_T)_x, q\right)$
      \STATE \label{algo_ressplitgcd_line2}$(S_2)_Q \gets (S_2)_Q \cup \{ q \}$
      \RETURN $S_1, S_2, \prs_i((S_T)_x, q, x)_{=}$
   \end{algorithmic}
\end{myalgorithm}
\ifthenelse{\equal{\articletype}{split} \or \equal{\splitarticle}{jsc}}{
\begin{proof}[Correctness]
 Property \ref{algo_ressplitgcd_cond1}) follows from Algorithm (\ref{algo_ressplit}) and line \ref{algo_ressplitgcd_line2}. Property \ref{algo_ressplitgcd_cond2}) was already shown in the correctness proof of the last algorithm.\qed
\end{proof}
}{}

\ifthenelse{\equal{\articletype}{split} \or \equal{\splitarticle}{jsc}}{
Remark that $i>0$ must always hold for this algorithm to make sense - this condition will always be true later. If $i=0$, then we will get $\widetilde{q} = \res_0((S_T)_x,q,x)$. As the latter has been added as an inequation to $S_1$ and we want to add $\widetilde{q}$ to $S_1$ as an equation, $i=0$ will always lead to a contradiction in $S_1$.
}

The following algorithm is similar, but instead of the gcd, it returns the first input polynomial divided by the gcd. It is used to assimilate an inequation into a system where there already is an equation with the same leader, or to calculate the least common multiple of two inequations.

\begin{myalgorithm}[\sf ResSplitDivide]\label{algo_ressplitdivide}
 \textit{Input:} A system $S$ with $(S_Q)_{<x}^==\emptyset$ and two polynomials $p$, $q$ with $\ld(p)=\ld(q)=x$ and $\pma(\ini(p))\neq0$ for all $\ma \in \sol(S_{<x})$.
                 Furthermore, if $\rank(p)\le\rank(q)$ then $\pma(\ini(q))\neq0$.\\
 \textit{Output:} Two systems $S_1$ and $S_2$ and a polynomial $\widetilde{p}$ such that:
   \renewcommand{\theenumi}{\alph{enumi}}
   \renewcommand{\labelenumi}{\theenumi)}
  \begin{enumerate}
   \item\label{algo_ressplitdivide_cond1} $S_2 = \widetilde{S_2} \cup \{ q \}$ where $\left(S_1,\widetilde{S_2}\right)$ is an $i$-th fibration split $p$ and $q^\prime$ w.r.t. $S$
   \item\label{algo_ressplitdivide_cond2} $\widetilde{p}$ is an $i$-th conditional quotient of $p$ by $q^\prime$ w.r.t. $S$
  \end{enumerate} where $i$ is the quasi fiber cardinality of $p$ and $q^\prime$ w.r.t. $S$, with $q^\prime=q$ for $\rank(p)>\rank(q)$ and $q^\prime=\prem(q,p,x)$ otherwise.\\
   \renewcommand{\theenumi}{\arabic{enumi}}
   \renewcommand{\labelenumi}{\theenumi)}
 \textit{Algorithm:} \begin{algorithmic}[1]
  \IF{$\rank(p)\le\rank(q)$}
   \RETURN $\mathsf{ResSplitDivide}(S, p, \prem(q,p,x))$
  \ELSE
   \STATE $(i, S_1, S_2) \gets \textrm{\sf ResSplit}\left(S, p, q\right)$
   \IF{$i>0$}
     \STATE $\widetilde{p} \gets \pquo(p, \prs_i(p, \prem(q,p,x), x),x)$
   \ELSE
     \STATE $\widetilde{p} \gets p$
   \ENDIF
   \STATE $(S_2)_Q \gets (S_2)_Q \cup \{ q \}$ \label{algo_ressplitdivide_line7}
   \RETURN $S_1, S_2, \widetilde{p}$
  \ENDIF
 \end{algorithmic}
\end{myalgorithm}
\ifthenelse{\equal{\articletype}{split} \or \equal{\splitarticle}{jsc}}{
 \begin{proof}[Correctness]
  According to Corollary (\ref{cor_ressplit}), we can without loss of generality assume $\rank(p)>\rank(q)$.

  Property \ref{algo_ressplitgcd_cond1}) follows from Algorithm (\ref{algo_ressplit}) and line \ref{algo_ressplitdivide_line7}.
  For all $\ma \in \sol(S_1)$, the following holds:
  If $i=0$, then $\deg_x(\gcd(\pma[<x](p), \pma[<x](q^\prime)))=0$ and thus $\pma[<x](p)$ shares no roots with $\pma[<x](q^\prime)$.
  Now let $i>0$. Formula (\ref{prempquo}) on page \pageref{prempquo} implies
    \[m\cdot p = \widetilde{p} \cdot \prs_i\left(p, q^\prime, x\right)+\prem\left(p,\prs_i\left(p, q^\prime, x\right),x\right)\]
  and due to \cite[Corr.~7.5.6]{mishra} and (\ref{lm_ressplit_proof1}), (\ref{lm_ressplit_proof2}) on page \pageref{lm_ressplit_proof1} there exist $k_1, k_2 \in F \setminus \{0\}$ such that:
   \[\begin{array}{rl} & \underbrace{\pma(m)}_{\neq 0}\cdot \pma[<x](p) \\
   = & \pma[<x](\widetilde{p}) \cdot \pma[<x]\left(\prs_i\left(p, q^\prime, x\right)\right)+\pma[<x]\left(\prem(p,\prs_i\left(p, q^\prime, x\right),x)\right)\\
   = & \pma[<x](\widetilde{p}) \cdot k_1 \prs_i\left(\pma[<x](p), \pma[<x](q^\prime), x\right)+k_2\prem(\pma[<x](p),\underbrace{\prs_i\left(\pma[<x](p), \pma[<x](q^\prime), x\right)}_{\mbox{\scriptsize divides } \pma[<x](p)},x)\\
   =& \pma[<x](\widetilde{p})\cdot k_1\gcd(\pma[<x](p), \pma[<x](q^\prime)) + 0 \end{array}\]
  and we obtain property \ref{algo_ressplitdivide_cond2}) from
   \[\pma[<x](\widetilde{p}) \sim \frac{\pma[<x](p)}{\gcd(\pma[<x](p), \pma[<x](q^\prime))}\]
  and $\deg_x(\pma[<x](\widetilde{p})) = \deg_x(\pma[<x](p))-\deg_x(\gcd(\pma[<x](p), \pma[<x](q^\prime)))=\deg_x(p)-i$.\qed
 \end{proof}
}{}

Applying the last algorithm to a polynomial $p$ and its partial derivative by its leader yields an algorithm to make polynomials square-free.
\ifthenelse{\equal{\articletype}{split} \or \equal{\splitarticle}{jsc}}{
We present it separately here to improve the readability of the \textsf{Decompose} algorithm.

\begin{myalgorithm}[\sf ResSplitSquareFree]\label{algo_ressplitsquarefree}
 \textit{Input:} A system $S$ with $(S_Q)_{<x}^==\emptyset$ and a polynomial $p$ with $\ld(p)=x$ and $\pma(\ini(p))\neq0$ for all $\ma \in \sol(S_{<x})$.\\
 \textit{Output:} Two systems $S_1$ and $S_2$ and a polynomial $r$ such that:
  \renewcommand{\theenumi}{\alph{enumi}}
  \renewcommand{\labelenumi}{\theenumi)}
  \begin{enumerate} \label{algo_ressplitsquarefree_cond1}
   \item $S_2 = \widetilde{S_2} \cup \{ p \}$ where $\left(S_1,\widetilde{S_2}\right)$ is an $i$-th square-free split of $p$ w.r.t. $S$
   \item \label{algo_ressplitsquarefree_cond2}$r$ is an $i$-th conditional square-free part of $p$ w.r.t. $S$
  \end{enumerate}
  \renewcommand{\theenumi}{\arabic{enumi}}
  \renewcommand{\labelenumi}{\theenumi)}
  where $i$ is the quasi fiber cardinality of $p$ and $\frac{\partial}{\partial x}p$ w.r.t. $S$.\\
 \textit{Algorithm:} \begin{algorithmic}[1]
  \STATE $(i, S_1, S_2) \gets \textrm{\sf ResSplit}\left(S, p, \frac{\partial}{\partial x}p\right)$
  \IF{$i>0$}
    \STATE $r \gets \pquo\left(p, \prs_i\left(p, \frac{\partial}{\partial x}p, x\right),x\right)$
  \ELSE
    \STATE $r \gets p$
  \ENDIF
  \STATE $(S_2)_Q \gets (S_2)_Q \cup \{ p \}$
  \RETURN $S_1, S_2, r$
 \end{algorithmic}
\end{myalgorithm}
 \begin{proof}[Correctness]
  Note that $\pma[<x](\frac{\partial}{\partial x}p)=\frac{\partial}{\partial x}\pma[<x](p)$ and therefore an $i$-th square-free split of $p$ is an $i$-th fibration split of $p$ and $\frac{\partial}{\partial x}p$. The rest follows from the proof of Algorithm (\ref{algo_ressplitdivide}).\qed
 \end{proof}
}{}

In the above \textsf{ResSplit}-based algorithms, we had the requirement that $(S_Q)_{<x}^==\emptyset$.
This ensures that all information contained in any equation of a smaller leader than $x$ will be respected by reduction modulo $S_T$ and thus avoids creating redundant systems.
It will also be necessary for \ifthenelse{\equal{\articletype}{split} \or \equal{\splitarticle}{jsc}}{proving} termination of the \textsf{Decompose} algorithm.
This motivates the definition of a selection strategy as follows.

\begin{mydefinition}[\sf Select]\label{Select}
 Let $\mathbb{P}_{\textrm{finite}}(M)$ be the set of all finite subsets of a set $M$. A \textbf{selection strategy} is a map \[\begin{array}{rcl} \textrm{\sf Select}: \mathbb{P}_{\textrm{finite}}(\{ p_=, p_{\neq} \mid p \in R \}) & \longrightarrow & \{ p_=, p_{\neq} \mid p \in R \}:\\ Q & \longmapsto & q \in Q\end{array}\] with the following properties: \begin{enumerate}
  \item\label{Select_Axiom1} If $\textrm{\sf Select}(Q) = q_=$ is an equation, then $Q_{<\ld(q)}^= = \emptyset$.
  \item\label{Select_Axiom2} If $\textrm{\sf Select}(Q) = q_{\neq}$ is an inequation, then $Q_{\leq\ld(q)}^= = \emptyset$.
 \end{enumerate}
\end{mydefinition}

The second property of \textsf{Select} could be weakened further, i.e., $Q_{<\ld(q)}^= = \emptyset$.
However, this would result in redundant calculations in the \textsf{Decompose} algorithm, thus we want all equations of the same leader to be treated first.

\ifthenelse{\equal{\articletype}{split} \or \equal{\splitarticle}{jsc}}{
We give a small counter-example to demonstrate what happens if the conditions of \textsf{Select} are violated.

\begin{myexample}
 Consider $R := F[a,x]$, the system $S$ with $S_T := \emptyset$ and $S_Q := \left\{ (x^2-a)_= \right\}$.
 To be able to insert $(x^2-a)_=$ into $S_T$ we need to apply the \textsf{ResSplitSquareFree} algorithm:
 We calculate $\res_0(x^2-a, 2x, x) = -4a$, $\res_1(x^2-a, 2x, x)=2$ and $\res_2(x^2-a, 2x, x)=1$ according to Definition (\ref{not_prs}).
 Thus, the quasi fiber cardinality is $0$ and we get the two new systems $S_1$, $S_2$ with \[(S_1)_T = \{ (x^2-a)_= \}, (S_1)_Q = \{ (-4a)_{\neq} \} \mbox{\ \ and\ \ } (S_2)_T = \emptyset, (S_2)_Q = \{ (x^2-a)_=, (-4a)_= \}\enspace.\]

 We now consider what happens with $S_2$: If we - violating the properties in Definition (\ref{Select}) - select $(x^2-a)_=$ as the next equation to be treated, $S_2$ will be split up into $S_{2,1}$, $S_{2,2}$ with \[(S_{2,1})_T = \{ (x^2-a)_= \}, (S_{2,1})_Q = \{ (-4a)_{\neq}, (-4a)_= \}\] and \[(S_{2,2})_T = \emptyset, (S_{2,2})_Q = \{ (x^2-a)_=, (-4a)_=, (-4a)_= \}\enspace.\] As $S_2 = S_{2,2}$, this will lead to an endless loop.
\end{myexample}
}{}

The following algorithm is trivial. However, it will be replaced with a more complicated algorithm in \crossref{\textsection}{section_differential}{} when the differential \textsc{Thomas} decomposition is treated.

\begin{myalgorithm}[\sf InsertEquation]
 \textit{Input:} A system $S$ and an equation $r_=$ with $\ld(r)=x$ satisfying $\pma(\ini(r))\neq0$ and $\pma[<x](r)$ square-free for all $\ma \in \sol(S_{<x})$. \\
 \textit{Output:} A system $S$ where $r_=$ is inserted into $S_T$. \\
 \textit{Algorithm:} \begin{algorithmic}[1]
  \IF{$(S_T)_x$ is not empty}
   \STATE $S_T \gets (S_T \setminus \{ (S_T)_x \})$
  \ENDIF
  \STATE $S_T \gets S_T \cup \{ r_= \}$
  \RETURN $S$
 \end{algorithmic}
\end{myalgorithm}

Now we present the main algorithm. It is based on all above algorithms and yields an algebraic \textsc{Thomas} decomposition. This algorithm \ifthenelse{\equal{\articletype}{split} \or \equal{\splitarticle}{jsc}}{also} forms the basis of the differential \textsc{Thomas} decomposition to be discussed in detail in \crossref{\textsection}{section_differential}{}.

The general structure is as follows:
In each iteration, a system $S$ is selected from a list $P$ of unfinished systems.
An equation or inequation $q$ is chosen from $S_Q$ according to the selection strategy and reduced modulo $S_T$.
The algorithm assimilates $q$ into $S_T$, potentially adding inequations of lower leader to $S_Q$ and adding new systems $S_i$ to $P$ that contain a new equation of lower leader in $(S_i)_Q$.
This process works differently depending on whether $q$ and $(S_T)_{\ld(q)}$ are equations or inequations, but it is based on the \textsf{InitSplit}, \textsf{ResSplitGcd}\ifthenelse{\equal{\articletype}{split} \or \equal{\splitarticle}{jsc}}{, \textsf{ResSplitDivide} and \textsf{ResSplitSquareFree}}{ and \textsf{ResSplitDivide}} methods in all cases.
As soon as the algorithm yields an equation $c=0$ for $c\in F\setminus\{0\}$ or an inequation $0\not=0$ in a system, this system is inconsistent and thus discarded.

\begin{myalgorithm}[\sf Decompose]\label{algo_decompose} The algorithm is printed on page \pageref{algo_decompose_float}.
 \begin{Decompose}
 \label{algo_decompose_float}
 \begin{minipage}[top]{\textwidth}
 {\textbf{Algorithm \ref{algo_decompose} (\textsf{Decompose})}}\vspace{-5pt}\\
 \noindent\begin{tabular*}{\textwidth}{c}\hline\end{tabular*} 
 \textit{Input:} A system $S^\prime$ with $({S^\prime})_T=\emptyset$. \\
 \textit{Output:} A {\sc Thomas} decomposition of $S^\prime$. \\
 \textit{Algorithm:}
 \begin{algorithmic}[1]
  \STATE $P \gets \{ S^\prime \}$; $\mathit{Result} \gets \emptyset$
  \WHILE{$|P|>0$}\label{main_loop}
   \STATE \label{choose_S} Choose $S \in P$; $P \gets P \setminus \{ S \}$
   \IF{$|S_Q|=0$}\label{SQempty}
    \STATE $\mathit{Result} \gets \mathit{Result} \cup \{ S \}$
   \ELSE
    \STATE \label{remove_from_Q} $q \gets \textrm{\sf Select}(S_Q)$; $S_Q \gets S_Q \setminus \{ q \}$
    \STATE \label{reduce_q} $q \gets \textrm{\sf Reduce}(q, S_T)$; $x \gets \ld(q)$
    \IF{$q \notin \{ 0_{\neq}, c_{=} \ |\ c \in F \setminus \{0\} \}$}\label{omit_sys}
    \IF{$x \neq 1$}
    \IF{$q$ is an equation}
     \IF{$(S_T)_x$ is an equation}
      \IF{$\textrm{\sf Reduce}(\res_0((S_T)_x, q, x), S_T)=0$}\label{reduceres}
       \STATE \label{resplitgcd} $(S, S_1, p) \gets \textrm{\sf ResSplitGCD}(S, q, x)$; $P \gets P \cup \{ S_1 \}$
       \STATE \label{q_1}$S \gets \textrm{\sf InsertEquation}(S, p_=)$
      \ELSE
       \STATE \label{insertres}$S_Q \gets S_Q \cup \{ q_=, \res_0((S_T)_x, q, x)_= \}$
      \ENDIF
     \ELSE
      \IF{$(S_T)_x$ is an inequation\footnote{Remember that $(S_T)_x$ might be empty, and thus neither an equation nor an inequation.}}\label{begin_eqineq}
       \STATE \label{removeineq} $S_Q \gets S_Q \cup \{ (S_T)_x \}$; $S_T \gets S_T \setminus \{ (S_T)_x \}$
      \ENDIF
      \STATE $(S, S_2) \gets \textrm{\sf InitSplit}(S, q)$; $P \gets P \cup \left\{ S_2 \right\}$
      \ifthenelse{\equal{\articletype}{split} \or \equal{\splitarticle}{jsc}}{
      \STATE $(S,S_3,p) \gets \textrm{\sf ResSplitSquareFree}\left(S, q\right)$; $P \gets P \cup \{ S_3 \}$
      }{
      \STATE $(S,S_3,p) \gets \textrm{\sf ResSplitDivide}\left(S, q, \frac{\partial}{\partial x}q\right)$; $P \gets P \cup \{ S_3 \}$
      }
      \STATE \label{q_2}$S \gets \textrm{\sf InsertEquation}(S, p_=)$
     \ENDIF
    \ELSIF{$q$ is an inequation}
     \IF{$(S_T)_x$ is an equation}
      \STATE \label{dividebyinequation} $(S, S_4, p) \gets \textrm{\sf ResSplitDivide}\left(S, (S_T)_x, q\right)$; $P \gets P \cup \{ S_4 \}$
      \STATE \label{q_3}$S \gets \textrm{\sf InsertEquation}(S, p_=)$
     \ELSE\label{begin_ineqineq}
      \STATE $(S, S_5) \gets \textrm{\sf InitSplit}(S, q)$; $P \gets P \cup \{ S_5 \}$
      \ifthenelse{\equal{\articletype}{split} \or \equal{\splitarticle}{jsc}}{
      \STATE $(S, S_6, p) \gets \textrm{\sf ResSplitSquareFree}\left(S, q\right)$; $P \gets P \cup \{ S_6 \}$
      }{
      \STATE $(S, S_6, p) \gets \textrm{\sf ResSplitDivide}\left(S, q, \frac{\partial}{\partial x}q\right)$; $P \gets P \cup \{ S_6 \}$
      }
      \IF{$(S_T)_x$ is an inequation}
       \STATE $(S, S_7, r) \gets \textrm{\sf ResSplitDivide}\left(S, (S_T)_x, p \right)$; $P \gets P \cup \{ S_7 \}$
       \STATE \label{q_4}$(S_T)_x \gets (r\cdot p)_{\neq}$
      \ELSIF{$(S_T)_x$ is empty}
       \STATE \label{q_5}$(S_T)_x \gets p_{\neq}$
      \ENDIF\label{end_ineqineq}
     \ENDIF
    \ENDIF
    \ENDIF
    \STATE \label{add_S}$P \gets P \cup \{ S \}$
    \ENDIF
   \ENDIF
  \ENDWHILE
  \RETURN $\mathit{Result}$
 \end{algorithmic}
 \end{minipage}
 \end{Decompose}
\end{myalgorithm}

\ifthenelse{\equal{\articletype}{split} \or \equal{\splitarticle}{jsc}}{
\begin{myexample}
 We demonstrate the algorithm with the simple example 
 \[S = (S_T, S_Q) := (\emptyset, \{ (x^2+x+1)_=, (x+a)_{\neq} \})\]
 with $a<x$.
 According to \textsf{Select}, $q := (x^2+x+1)_=$ is chosen in line \ref{remove_from_Q}.

 As $\ini(q)=1$, the system $\widetilde{S} := S \cup \{ 1_= \}$ is added to the set $P$ of systems to be processed.
 However, as $1_= \in \widetilde{S}$, this system does not have any solution.
 In a later iteration of the loop, this system will be chosen again in line \ref{choose_S} and not put back into $P$ as line \ref{add_S} is skipped due to the condition in line \ref{omit_sys} applied to $1_=$.
 For readability, we will directly omit such systems for the rest of the example.
 
 The same is true for the system where $\res_0(q, \frac{\partial}{\partial x}q, x)=1$ is added as an equation.
 Therefore, in line \ref{q_2}, the original system $S$ is replaced with $\left(\{(x^2+x+1)_=\}, \{(x+a)_{\neq}\}\right)$.

 Now, $q := (x+a)_{\neq}$ is chosen and $\mathsf{ResSplitDivide}(S, (S_T)_x, q)$ is computed in line \ref{dividebyinequation}:
 By Definition (\ref{not_prs}), $\res_0((S_T)_x, q, x) = \prem((S_T)_x, q, x) = a^2-a+1$, $\res_1((S_T)_x, q, x)=\ini(q)=1$, and $\res_2((S_T)_x, q, x)=1$.
 $S_T$ contains no equation of leader $a$, so no reduction is possible.
 Thereby, $S$ is decomposed into
  \[S := (\underbrace{\{ (x^2+x+1)_=, (a^2-a+1)_{\neq} \}}_{= S_T}, \underbrace{\{\}}_{=S_Q})\]
 which is already simple and
  \[S_1 := ( \underbrace{\{ (x^2+x+1)_= \}}_{=(S_1)_T}, \underbrace{\{ (x+a)_{\neq}, (a^2-a+1)_= \}}_{=(S_1)_Q} )\enspace.\]
 Omitting trivial cases again, $S_1$ is replaced with \[S_1 := \left( \{ (x^2+x+1)_=, (a^2-a+1)_= \}, \{ (x+a)_{\neq} \} \right)\enspace.\]
 For $S_1$, the procedure $\mathsf{ResSplitDivide}(S_1, ((S_1)_T)_x, q)$ needs to be performed again.
 This time $\mathsf{Reduce}(a^2-a+1, (S_1)_T)=0$ holds and - again omitting trivial cases - $S_1$ is replaced with \[S_1 := ( \{ \underbrace{(x-a+1)_=}_{\pquo(x^2+x+1, x+a, x)}, (a^2-a+1)_= \}, \{ 1_{\neq} \} )\enspace.\]
 Finally, the {\sc Thomas} decomposition of $S$ is:
 \[\left(\{ (x^2+x+1)_=, (a^2-a+1)_{\neq} \}, \{ (x-a+1)_=, (a^2-a+1)_= \} \right)\enspace.\]
\end{myexample}
}{}

\ifthenelse{\equal{\articletype}{split} \or \equal{\splitarticle}{jsc}}{
 \begin{proof}[Correctness]
  
  The correctness of the \textsf{Decompose} algorithm is proved by verifying two loop invariants:
  \begin{enumerate}
   \item $P \cup \mathit{Result}$ is a disjoint decomposition of the input $S^\prime$.
   \item For all systems $S \in P \cup \mathit{Result}$ it holds that $S_T$ is triangular and for all $p \in S_T$ with $\ld(p)=x$, $\pma[<x](p)$ is square-free and $\pma(\ini(p))\neq0$ for all solutions $\ma \in \sol((S_T)_{<x} \cup (S_Q)_{<x})$.
  \end{enumerate}

  We prove the first loop invariant.
  At the beginning $P\cup\mathit{Result}=\{S^\prime\}$ holds, which is clearly a disjoint decomposition of $S^\prime$.
  Now assume that $P \cup \mathit{Result}$ is a disjoint decomposition of $S^\prime$ at the beginning of the loop.
  It suffices to show that all systems we add to $P$ or $\mathit{Result}$ until the end of the loop add up to a disjoint decomposition of the system $S$ chosen in line \ref{choose_S}.
  If in line \ref{SQempty} $S_Q=\emptyset$ holds, the algorithm just moves $S$ from $P$ to $\mathit{Result}$, the statement is obvious.

 In line \ref{insertres}, $S$ is not decomposed, but rather a new equation is added.
 However, if $p,q \in R$ with $\ld(p)=x$, $\pma[<x](p)=0$ and $\pma[<x](q)=0$ for each $\ma \in F[\ y\mid y<x\ ]$ imply $\pma(\res_0(p,q,x))=0$ (cf.\  \cite[Lemma 7.2.3]{mishra}).
 Therefore, we do not loose any solutions by adding the equation $\res_0((S_T)_x,q,x)_=$.

  Note now that if $(S, S_i)$ is the output of any of the {\sf ResSplitGcd}, {\sf InitSplit}, {\sf ResSplitSquareFree} and {\sf ResSplitDivide} algorithms, then $(S \cup \{q\}, S_i)$ is a disjoint decomposition of $S \cup \{q\}$, where $S$ is the input system of the respective algorithm.
  It remains to be shown that the actions performed in lines \ref{q_1}, \ref{q_2}, \ref{q_3}, \ref{q_4} and \ref{q_5} are equivalent to putting $q$ back into the system $S$.

  For any $\ma \in \sol(S_{<x})$, the following holds: For line \ref{q_1}, Algorithm (\ref{algo_ressplitgcd}) guarantees
   \[\pma[<x](p)=0 \Longleftrightarrow \pma[<x]((S_T)_x)=0 \textrm{ and } \pma[<x](q)=0\] where $p_=$ is a minimal conditional gcd.
  For line \ref{q_3}, Algorithm (\ref{algo_ressplitdivide}) ensures that
   \[\pma[<x](p)=0 \Longleftrightarrow \pma[<x]((S_T)_x)=0 \textrm{ and } \pma[<x](q)\neq0\enspace.\]
  For line \ref{q_2}, (\ref{algo_ressplitsquarefree}) guarantees that $p$ has the same solutions as $q$ because
   \[\pma[<x](p) \sim \frac{\pma[<x](q)}{\gcd(\pma[<x](q), \pma[<x](\frac{\partial}{\partial x}q))} = \frac{\pma[<x](q)}{\gcd(\pma[<x](q), \frac{\partial}{\partial x}\pma[<x](q))}\enspace.\]
  In lines \ref{q_4} and \ref{q_5}, the roots of $p$ and $q$ are identical with the same reasoning as before. Furthermore, in line \ref{q_4},
   \[\pma[<x](r) \sim \frac{\pma[<x]((S_T)_x)}{\gcd(\pma[<x]((S_T)_x), \pma[<x](p))}
  \Longrightarrow
   \pma[<x](r\cdot p) \sim \lcm(\pma[<x]((S_T)_x), \pma[<x](p))\enspace.\]
  This concludes the proof of the first loop invariant.

  Now, we prove the second loop invariant:
  At the start, the only system in $P$ is $S^\prime$ with $S^\prime_T=\emptyset$.
  So assume that the second loop invariant holds at the beginning of a loop.
  
  One easily checks that all steps in the algorithm allow only one polynomial $(S_T)_x$ in $S_T$ for each leader $x$, thus triangularity obviously holds.

  We show that all polynomials added to $S_T$ have non-zero initial and are square-free.
  For $\sol(S_{<x})=\emptyset$, the statement is trivially true.
  Otherwise, let again $\ma \in \sol(S_{<x})$.

  For the equation $p_=$ added as conditional gcd of $(S_T)_x$ and $q$ in line \ref{q_1}, it holds that $\pma[<x](p)$ is a divisor of $\pma[<x]((S_T)_x)$. As $\pma[<x]((S_T)_x)$ is square-free by assumption, so is $\pma[<x](p)$.
  As $\textrm{\sf Reduce}(res_0((S_T)_x, q, x))=0$, the quasi fiber cardinality of $(S_T)_x$ and $q$ is positive, thus the inequation added to $S$ in {\sf ResSplit} is by Definition (\ref{not_prs}) the initial of $p_=$.

  The equation $p_=$ inserted into $S_T$ in line \ref{q_2} and the inequation $p_{\neq}$ inserted in line \ref{q_5} are square-free due to Algorithm (\ref{algo_ressplitsquarefree}) and their initials are non-zero as $p$ is either identical to $q$, or it is a pseudo quotient of $q$ by $\prs_i\left(q,\frac{\partial}{\partial x}q, x\right)$ for some $i>0$.
  On the one hand, if $p$ equals $q$, the call of \textsf{InitSplit} for $q$ ensures a non-zero initial for $p$.
  On the other hand, the polynomial $\prs_i\left(q,\frac{\partial}{\partial x}q, x\right)$ has initial $\res_i\left(q,\frac{\partial}{\partial x}q, x\right)$, which is added as an inequation by {\sf ResSplit}.
  This implies that the initial of the pseudo-quotient is also non-zero, since again the initial of $q$ was ensured to be non-zero by \textsf{InitSplit}.

  The equation $p_=$ that replaces the old equation $(S_T)_x$ in line \ref{q_3} is the quotient of $(S_T)_x$ by an inequation.
  It is square-free, because $\pma[<x](p)$ is a divisor of $\pma[<x]((S_T)_x)$ being square-free by assumption.
  Also $p$ is either identical to $(S_T)_x$ or a pseudo quotient of $(S_T)_x$ by $\prs_i\left((S_T)_x, q, x\right)$ for some $i>0$, where the initial of the pseudo quotient is non-zero for the same argument as above.

  Finally, consider the inequation $(r\cdot p)_{\neq}$ added in line \ref{q_4} as a ``least common multiple'' of $\left((S_T)_x\right)_{\not=}$ and $p_{\not=}$.
  The inequation $\pma[<x](p)$ is square-free and has non-vanishing initial for the same reasons as before (see the above comments about $p_{\neq}$ from line \ref{q_5}).
  Due to $\pma[<x](r) \sim \frac{\pma[<x]((S_T)_x)}{\gcd(\pma[<x]((S_T)_x), \pma[<x](p))}$, the polynomials $\pma[<x](r)$ and $\pma[<y](p)$ have no common divisors.
  As $\pma[<x](r)$ divides $\pma[<x]((S_T)_x)$, using the same arguments as before, $\pma[<x](r)$ is square-free and has a non-vanishing initial.
  This completes the proof of the second loop invariant.

  It is obvious that a system $S$ with $S_Q=\emptyset$ for which these loop invariants hold is simple. Thus the algorithm returns the correct result if it terminates.\qed
 \end{proof}
}{}

\ifthenelse{\equal{\articletype}{split} \or \equal{\splitarticle}{jsc}}{
We now aim at showing termination and remark that the system $S$ chosen from $P$ is treated in one of three ways:
It might be discarded, put into $\mathit{Result}$, or replaced in $P$ by at least one new system.
To show that $P$ is empty after finitely many iterations we define orders on the systems and show that these orders are well-founded.
Afterwards we show in detail that the algorithm produces descending chains of systems, proving termination.

\begin{mydefinition}
For transitive and asymmetric\footnote{A relation $\prec$ is asymmetric, if $S\prec S'$ implies $S'\not\prec S$ for all $S,S'$. Asymmetry implies irreflexivity.} partial orders $<_i$ for $i=1,\ldots,m$, we define the \textbf{composite order} $\mbox{``}<\mbox{''}:=[<_1,\ldots,<_m]$ as follows:
$a<b$ if and only if there exists $i \in \{ 1, \ldots, m \}$ such that $a<_ib$ and neither $a<_jb$ nor $b<_ja$ for $j<i$.
The composite order is clearly transitive and asymmetric.
An order $<$ is called \textbf{well-founded}, if each $<$-descending chain becomes stationary.
\end{mydefinition}

\noindent The following trivial statement will be used repeatedly: 
  
\begin{myremark}\label{composite_well_founded}
If each $<_i$ is well-founded, then so is the composite ordering $<$.  
\end{myremark}

\noindent Now we define the orderings and show their well-foundedness:
  
\begin{defrem}\label{def_orders}
   
  Define the $\prec$ as the composite order $[\prec_1,\prec_2,\prec_3,\prec_4]$ of the four orders defined below. It is well-founded since the $\prec_i$ are.

   \begin{enumerate}

   \item \label{order_1}For $i=1,\ldots,n$ define $\prec_{1,x_i}$ by $S\prec_{1,x_i} S'$ if and only if $\rank\left((S_T)^=_{x_i}\right) < \rank\left((S^\prime_T)^=_{x_i}\right)$, with $\rank\left((S_T)^=_{x_i}\right):=\infty$ if $(S_T)^=_{x_i}$ is empty.
   Define the composite order $\prec_1$ as $[\prec_{1,x_1},\ldots,\prec_{1,x_n}]$.\\
   Since degrees can only decrease finitely many times, the orders $\prec_{1,x_i}$ are clearly well-founded and thus $\prec_1$ is, too.

   \item \label{order_2}Define the map $\mu$ from the set of all systems over $R$ to $\{1,x_1,\ldots,x_n,x_\infty\}$, where $\mu(S)$ is minimal such that there exists an equation $p \in (S_Q)^=_{\mu(S)}$ with $\mathsf{Reduce}(p, S_T)\neq0$, or $\mu(S)=x_\infty$ if no such equation exists. Then, $S \prec_2 S^\prime$ if and only if $\mu(S) < \mu(S^\prime)$ with $1<x_\infty$ and $x_i<x_\infty$ for $i\in\{1,\ldots,n\}$.

   The ordering $\prec_2$ is well-founded since the ranking $<$ is well-founded on $\{1,x_1, \ldots, x_n,x_\infty\}$, which is a finite set.

  \item \label{order_3}$S \prec_3 S'$ if and only if exists $p_{\not=}\in R^{\not=}$ and a finite set $L\subset R^{\not=}$ with $\ld(q)<\ld(p)\ \forall\ q\in L$ such that $S_Q \uplus \{p_{\not=}\} = S'_Q \uplus L$ holds.\\
  We show well-foundedness by induction on the highest appearing leader $x$ in $(S_Q)^{\not=}$:\\
  For $x=1$ we can only make a system $S$ $\prec_3$-smaller by removing one of the finitely many inequations in $(S_Q)^{\not=}$.\\
  Now assume that the statement is true for all indeterminates $y<x$.\\
  By the induction hypothesis we can only $\prec_3$-decrease $S$ finitely many times without changing $(S_Q)_x^{\not=}$.
  To further $\prec_3$-decrease $S$, we have to remove an inequation in $(S_Q)^{\not=}_x$ (and possibly replace it by one or more lower ranking inequations).
  As $(S_Q)^{\neq}_x$ is finite, these processes can only be repeated finitely many times until $(S_Q)^{\neq}_x=\emptyset$. Now, the highest appearing leader in $(S_Q)^{\neq}$ is smaller than $x$ and by the induction hypothesis, the statement is proven.

   \item \label{order_4}$S \prec_4 S^\prime$ if and only if $|S_Q| < |S^\prime_Q|$, which is again well-founded.
  \end{enumerate}
  
\end{defrem}

 \begin{proof}[Termination]
  
  We denote the system chosen from $P$ in line \ref{choose_S} by $\widehat{S}$ and the system added to $P$ in line \ref{add_S} by $S$.
  We prove that the systems $S, S_1, \ldots, S_7$ generated from $\widehat{S}$ are $\prec$-smaller than $\widehat{S}$.
  For $i=1,\ldots,4$ we will use the notation $S \osim_i S^\prime$ if neither $S \prec_i S^\prime$ nor $S^\prime \prec_i S$ holds.

  For $j=1,\ldots,7$, $((S_j)_T)^= = (\widehat{S}_T)^=$ and thus $S_j \osim_1 \widehat{S}$.
  The properties of {\sf Select} in Definition (\ref{Select}) directly require, that there is no equation in $(\widehat{S}_Q)^=$ with a leader smaller than $x$.
  However, the equation added to the system $S_j$ returned from {\sf InitSplit} (\ref{algo_initsplit}) is the initial of $q$, which has leader smaller $x$ and does not reduce to $0$, as ensured by {\sf Reduce} (\ref{algo_reduce}).
  Furthermore, the equations added in one the subalgorithms based on {\sf ResSplit} (\ref{algo_ressplit}) have a leader smaller than $x$ and do not reduce to $0$ due to line \ref{ressplit_cond_nontrivial_split} in algorithm {\sf ResSplit}.
  In each case $S_j \prec_2 \widehat{S}$ is proven.

  As soon as we show $S \prec \widehat{S}$, the proof is finished.
  If $q$ is reduced to $0_=$, then an equation that reduces to zero is omitted from $S_Q$, implying $\left|S_Q\right| < \left|\widehat{S}_Q\right|$ and thus $S\prec_4\widehat{S}$.
  This is equivalent to $S\prec\widehat{S}$, since $S\osim_i\widehat{S}, i=1,2,3$.
  If $q$ is reduced to $c_{\neq}$ for some $c \in F \setminus \{0\}$, then $S \prec_3 \widehat{S}$ and $S \osim_i \widehat{S}, i=1,2$, since the only change was the removal of an inequation from $S_Q$.

  Lines \ref{resplitgcd}-\ref{q_1} set $(S_T)_x$ to $p_=$ of smaller degree than $(\widehat{S}_T)_x$ and \ref{begin_eqineq}-\ref{q_2} add $(S_T)_x$ as a new equation.
  In both cases we get $S \prec_1 \widehat{S}$.
  In line \ref{insertres}, $S_T=\widehat{S}_T$ implies $S \osim_1 \widehat{S}$ and $q$ is chosen according to {\sf Select} (cf.\ \ref{Select}.\ref{Select_Axiom1}), which implies $(\widehat{S}_Q)^=_{<x}=\emptyset$ and $(S_Q)^=_{<x}=\{\res_0((S_T)_x, q, x)_=\}$.
  Because $\mbox{\sf Reduce}(\res_0((S_T)_x, q, x), S) \neq 0$ was ensured by line \ref{reduceres}, $S \prec_2 \widehat{S}$ follows.
  
  In lines \ref{dividebyinequation}-\ref{q_3}, if the degree of $(S_T)_x$ is smaller than the degree of $(\widehat{S}_T)_x$, then $S \prec_1 \widehat{S}$.
  In case the degree doesn't change we have $S \osim_1 \widehat{S}$ by definition and $(S_Q)^==(\widehat S_Q)^=$ guarantees $S \osim_2 \widehat{S}$.
  However, $q$ is removed from $S_Q$ and replaced by an inequation of smaller leader, thus implying $S \prec_3 \widehat{S}$.
  In \ref{begin_ineqineq}-\ref{end_ineqineq}, obviously $S \osim_i \widehat{S}, i=1,2$.
  As before, $q$ is removed from $S_Q$ and replaced by an inequation of smaller leader, which once more implies $S \prec_3 \widehat{S}$.\qed
 \end{proof}
}{}

 In the next section, we consider an extension of this algorithm to partial differential systems.
 Both algorithms have been implemented, and their implementation aspects are considered in \textsection\ref{implementation}.

%% file: curve
\begin{tikzpicture}[domain=-3:3, smooth, x=1.2cm, y=1.2cm]

\draw[color=gray, line width=1pt, style=dashed] plot[] coordinates{
(0,-0.75)
(0,0.75)
}
node[left]{\scriptsize\color{black}$x_=$};

\draw[color=gray, line width=1pt, style=dashed] plot[] coordinates{
(-1,-0.75)
(-1,0.75)
}
node[left]{\scriptsize\color{black}$(x+1)_=$};

\draw[color=black, line width=1pt] plot[] coordinates{
(0.589436386046560012, 0.743896947953439969)
(0.575757575999999993, 0.722126569068546553)
(0.565117125970401091, 0.707610147029598902)
(0.557317740107367898, 0.696969696999999999)
(0.540371889494211222, 0.671749322921924907)
(0.516595852484581175, 0.636363636999999982)
(0.516048157718607348, 0.635466994281392750)
(0.515151515000000004, 0.633999081467944814)
(0.490706749754113769, 0.600202341649224924)
(0.473026152128687594, 0.575757575999999993)
(0.465739372041181543, 0.564563658774118848)
(0.454545454999999987, 0.547367616175383875)
(0.440375618576493211, 0.529321351423506892)
(0.429249534202270744, 0.515151515000000004)
(0.414735128721764323, 0.494355780621365315)
(0.393939393999999998, 0.464560376595814284)
(0.389302937533814342, 0.459181911466185699)
(0.385306122820929964, 0.454545454999999987)
(0.362959516351743672, 0.424919272159424266)
(0.339591836842798145, 0.393939393999999998)
(0.337134358285293867, 0.390138369714706112)
(0.333333333999999981, 0.384259260180555562)
(0.310325157267228924, 0.356341510353136104)
(0.291363163872195630, 0.333333333999999981)
(0.283775447664895264, 0.322285159335104709)
(0.272727272999999992, 0.306198347428698336)
(0.256730959360571553, 0.288723586639428431)
(0.242088009271773485, 0.272727272999999992)
(0.229419899617189338, 0.255428585382810658)
(0.212121212000000003, 0.231806703236502099)
(0.202059315309927295, 0.222183108856094019)
(0.191539153853590799, 0.212121212000000003)
(0.173948138612201764, 0.189688225387798226)
(0.151515151999999986, 0.161080502556703364)
(0.146172586622852235, 0.156857717288995391)
(0.139413302754830604, 0.151515151999999986)
(0.117220734469718771, 0.125203508530281227)
(0.0909090909999999974, 0.0940082645556559876)
(0.0889086464525618569, 0.0929095355474381379)
(0.0852664577441141491, 0.0909090909999999974)
(0.0590744993961474785, 0.0621376216038525206)
(0.0303030300000000016, 0.0305325984049012907)
(0.0300751876770309241, 0.0305308723267284747)
(0., 0.)
};

\draw[color=black, line width=1pt] plot[] coordinates{
(0., 0.)
(0.0303030300000000016, -0.0305325984049012872)
(0.0852664577441141491, -0.0909090909999999974)
(0.0909090909999999974, -0.0940082645684400670)
(0.139413301700517522, -0.151515150999999987)
(0.151515151999999986, -0.161080502416341448)
(0.191539153853590799, -0.212121212000000003)
(0.212121212000000003, -0.231806703236502099)
(0.242088009271773485, -0.272727272999999992)
(0.272727272999999992, -0.306198347483925604)
(0.291363163116732604, -0.333333333000000009)
(0.333333333999999981, -0.384259260167245376)
(0.339591836842798145, -0.393939393999999998)
(0.358803984955309574, -0.419410044955309591)
(0.385306122012766705, -0.454545454000000015)
(0.393939393999999998, -0.464560376543642239)
(0.429249534202270744, -0.515151515000000004)
(0.454545454999999987, -0.547367616175383875)
(0.473026152128687594, -0.575757575999999993)
(0.515151515000000004, -0.633999081515994045)
(0.516595851843211551, -0.636363636000000010)
(0.519553788506024739, -0.640765909506024745)
(0.557317740107367898, -0.696969696999999999)
(0.575757575999999993, -0.722126569085841830)
} node[right]{$p_=$};

\draw[color=black, line width=1pt] plot[] coordinates{
(-0.939393939000000011, -0.229625804093301689)
(-0.907830033866470654, -0.272727272999999992)
(-0.878787878999999994, -0.304453627015615180)
(-0.841044104850747876, -0.333333333000000009)
(-0.818181818000000005, -0.347910927527046054)
(-0.794078571984453530, -0.357436579015546485)
(-0.757575757000000016, -0.371862565194621464)
(-0.710411690152531472, -0.380497400625675775)
(-0.696969696999999999, -0.382958371582926782)
(-0.647281921348538569, -0.383021108651461328)
(-0.636363636000000010, -0.383034894360465450)
(-0.593154375818473079, -0.376542593894479749)
(-0.575757575999999993, -0.373928680786649648)
(-0.543829556688303017, -0.365261352311696985)
(-0.515151515000000004, -0.357476277863315928)
(-0.497430280223396193, -0.351054567776603821)
(-0.454545454000000015, -0.335514232997603856)
(-0.453012825093495997, -0.334865961931792355)
(-0.449389415801167369, -0.333333333000000009)
(-0.412608496824779181, -0.314664230175220827)
(-0.393939393999999998, -0.305188246160087218)
(-0.372914622151552355, -0.293752044501538967)
(-0.334261838895042740, -0.272727272999999992)
(-0.333709131890942534, -0.272351474109057412)
(-0.333333333000000009, -0.272095959332899262)
(-0.297632469094419039, -0.247822076494645188)
(-0.272727272999999992, -0.230888429912351506)
(-0.261245128088579137, -0.223603356911420859)
(-0.243147896701875699, -0.212121212000000003)
(-0.226340007329028481, -0.197902416670971470)
(-0.212121212000000003, -0.185873890297195421)
(-0.192200185999680584, -0.171436177000319406)
(-0.164713334810353385, -0.151515150999999987)
(-0.158389704148128613, -0.144640597965301487)
(-0.151515150999999987, -0.137167125409578733)
(-0.126020117426885320, -0.116404124573114665)
(-0.0947145878450189188, -0.0909090909999999974)
(-0.0932669636705811128, -0.0885512182905140027)
(-0.0909090909999999974, -0.0847107438631198306)
(-0.0622962693675179757, -0.0589158516324820233)
(-0.0305576773062636829, -0.0303030300000000016)
(-0.0305576773062636829, -0.0300483826979379975)
(0,0)
};

\draw[color=black, line width=1pt] plot[] coordinates{
(0.0,0.0)
(-0.151515150999999987, 0.137167125218764341)
(-0.164713335960179175, 0.151515151999999986)
(-0.212121212000000003, 0.185873890369375572)
(-0.243147896701875699, 0.212121212000000003)
(-0.272727272999999992, 0.230888429912351506)
(-0.331360944841112115, 0.270754884808567742)
(-0.333333333000000009, 0.272095959332899262)
(-0.334261838895042740, 0.272727272999999992)
(-0.393939393999999998, 0.305188246106526617)
(-0.449389417968331484, 0.333333333999999981)
(-0.454545454000000015, 0.335514233077938484)
(-0.457965794820946959, 0.336753674764511290)
(-0.515151515000000004, 0.357476277913452711)
(-0.548290530602897341, 0.366472349056103597)
(-0.575757575999999993, 0.373928680814164416)
(-0.623531000483281028, 0.381106758483281016)
(-0.636363636000000010, 0.383034894375459123)
(-0.686002521855570979, 0.382972219036529449)
(-0.696969696999999999, 0.382958371598025704)
(-0.738915299789757896, 0.375278936789757989)
(-0.757575757000000016, 0.371862565224977015)
(-0.785191278746980581, 0.360948855291324533)
(-0.818181818000000005, 0.347910927590335206)
(-0.827083478943560602, 0.342234994796683123)
(-0.841044103413124189, 0.333333333999999981)
(-0.862426843243773833, 0.316972298243773820)
(-0.878787878999999994, 0.304453626963266721)
(-0.893950365222353494, 0.287889759472534534)
(-0.907830033866470654, 0.272727272999999992)
(-0.921173302142754391, 0.254506635842113860)
(-0.939393939000000011, 0.229625804093301689)
(-0.945601074414459619, 0.218328347414459611)
(-0.949011446425497618, 0.212121212000000003)
(-0.963890935744479593, 0.176012148340279034)
(-0.973985431653519473, 0.151515151999999986)
(-0.979591836601099497, 0.131106988601099483)
(-0.990634755436474057, 0.0909090909999999974)
(-0.991765782226405435, 0.0826748732264053765)
(-0.998959417293612151, 0.0303030300000000016)
(-0.998959417293612151, 0.0292624473107818001)
(-0.998959417293612151, -0.0303030300000000016)
(-0.998793727405677179, -0.0315093025943227953)
(-0.990634755436474057, -0.0909090909999999974)
(-0.987087517982503382, -0.103821572804440695)
(-0.973985431996911677, -0.151515150999999987)
(-0.955752212681415969, -0.195762938318584073)
(-0.949011446425497618, -0.212121212000000003)
(-0.939393939000000011, -0.229625804093301689)
};
\end{tikzpicture}

%% file: differential.tex
\section{The Differential Thomas Decomposition}\label{section_differential}

The differential \textsc{Thomas} decomposition is concerned with manipulations of polynomial differential equations.
The basic idea for a construction of this decomposition is twofold.
On the one hand a combinatorial calculus developed by {\sc Janet} takes care of finding unique reductors and all differential consequences by completing systems to involution.
On the other hand the algebraic {\sc Thomas} decomposition makes the necessary splits for regularity of initials during the computation and ensures disjointness.

\ifthenelse{\equal{\articletype}{combined}}{%
We start by giving the basic definitions from differential algebra needed for the algorithms. Afterwards we summarize the {\sc Janet} division and its combinatorics.
The combinatorics give us a new algorithm {\sf InsertEquation} to add equations into systems.
Afterwards we review the differential implications of the algebraic decomposition algorithm and present the algorithm {\sf Reduce} utilized for differential reduction.
Replacing the insertion and reduction methods from the algebraic case with these differential methods yields the differential decomposition algorithm.
}%

\subsection{Preliminaries from Differential Algebra}\label{differential_preliminaries}

Let $\Delta=\{\partial_1,\ldots,\partial_n\}$ be the set of derivations $(n>0)$ and $F$ be a computable $\Delta$\textbf{-differential field} of characteristic zero.
This means any $\partial_j\in\Delta$ is a linear operator $\partial_j:F\to F$ fulfilling the \textsc{Leibniz} rule.
For a \textbf{differential indeterminate} $u$ consider the $\Delta$\textbf{-differential polynomial ring} $F\{u\}:=F\left\lbrack\ u_{\mathbf{i}}\mid\mathbf{i}\in\Z_{\ge0}^n\ \right\rbrack$, a polynomial ring infinitely generated by the algebraically independent set $\langle u\rangle_\Delta:=\{u_{\mathbf{i}}\mid\mathbf{i}\in\Z_{\ge0}^n\}$.
The operation of $\partial_j\in\Delta$ on $\langle u\rangle_\Delta$ by $\partial_j u_{\mathbf{i}}=u_{\mathbf{i}+e_j}$ is extended linearly and by the \textsc{Leibniz} rule to $F\{u\}$.
Let $U=\{u^{(1)},\ldots,u^{(m)}\}$ be the set of differential indeterminates.
The multivariate $\Delta$-differential polynomial ring is given by $F\{U\}:=F\{u^{(1)}\}\ldots\{u^{(m)}\}$.
The elements of $\langle U\rangle_\Delta:=\left\{u^{(j)}_{\mathbf{i}}\mid\mathbf{i}\in\Z_{\ge0}^n, j\in\{1,\ldots,m\}\right\}$ are called \textbf{differential variables}.

We remark, that the algebraic closure $\overline{F}$ of $F$ is a differential field with a differential structure uniquely defined by the differential structure of $F$ (cf.\ \cite[\textsection II.2, Lemma 1]{Kol}).
Let
\[
  E:=
     \bigoplus_{j=1}^m\overline{F}[[z_1,\ldots,z_n]]
     \cong\overline{F}^{\langle U\rangle_\Delta}
\]
with indeterminates $z_1,\ldots,z_n$, where $\overline{F}[[z_1,\ldots,z_n]]$ denotes the ring of formal power series.
The isomorphism maps coefficients of the power series to function values of differential variables, i.e.,
\[
  \alpha:
  \bigoplus_{i=1}^m\overline{F}[[z_1,\ldots,z_n]]
  \to\overline{F}^{\langle U\rangle_\Delta}:
  \left(
    \sum_{\mathbf{i}\in\Z_{\ge0}^n}a_{\mathbf{i}}^{(1)}\frac{z^{\mathbf{i}}}{\mathbf{i}!},
    \ldots,
    \sum_{\mathbf{i}\in\Z_{\ge0}^n}a_{\mathbf{i}}^{(m)}\frac{z^{\mathbf{i}}}{\mathbf{i}!}
  \right)
  \mapsto\left(u^{(j)}_{\mathbf{i}}\mapsto a_{\mathbf{i}}^{(j)}\right)
\]
where $\mathbf{i}!:=i_1!\cdot\ldots\cdot i_n!$ defines the factorial of a multi-index.

In the formulation of the algorithm the direct sum of formal power series $E$ suffices to give a notion of solutions coherent to the algebraic case: 
For $e\in E$ we define the $F$-algebra homomorphism
\[
  \phi_{e}: F\{U\}\to\overline{F}: u^{(j)}_{\mathbf{i}}\mapsto\alpha(e)(u^{(j)})
\]
evaluating all differential variables of a differential polynomial at the power series $e$.
A \textbf{differential equation} or \textbf{inequation} for $m$ functions $U=\{u^{(1)},\ldots,u^{(m)}\}$ in $n$ indeterminates is an element $p\in F\{U\}$ written $p_=$ or $p_{\not=}$, respectively.
A \textbf{solution} of $p_=$ or $p_{\not=}$ is an $e\in E$ with $\phi_e(p)=0$ or $\phi_e(p)\not=0$, respectively.
More generally $e\in E$ is called a solution of a set $P$ of equations and inequations, if it is a solution of each element in $P$.
The set of solutions of $P$ is denoted by $\sol_E(P)=\sol(P)\subseteq E$.
Since we substitute elements of $\overline{F}$ algebraically for the differential indeterminates, \crossref{Remark (}{exist_sol}{)}, which guarantees the continuation of solutions from lower ranking variables to higher ranking ones, also holds here.

Any differential $F$-algebra $R$ with a differential embedding of $E\hookrightarrow R$ might be chosen as universal set of solutions, for example a universal differential field containing $F$:
Clearly $\overline{F}[[z_1,\ldots,z_n]]$ embeds into its field of quotients $\overline{F}((z_1,\ldots,z_n))$, and thus $\overline{F}[[z_1,\ldots,z_n]]$ also embeds into a universal differential field containing $F$, since $\overline{F}((z_1,\ldots,z_n))$ is a finitely generated differential field extension of $\overline{F}$ (cf.\ \cite[\textsection II.2 and \textsection III.7]{Kol}).  
We denote the set of solutions in $R$ by $\sol_R(P)\subseteq R$.

A finite set of equations and inequations is called a \textbf{(differential) system} over $F\{U\}$.
We will be using the same notation for systems as in the algebraic {\sc Thomas} decomposition introduced in \crossref{\textsection}{algebraic_definition_notation}{} and \crossref{\textsection}{algebraic_algorithms}{}, in particular a system $S$ is represented by a pair $(S_T,S_Q)$.
However, the candidate simple system $S_T$ will also reflect a differential structure using combinatorial methods.
We will elaborate on the combinatorics in the next section.

\subsection{The Combinatorics of Janet Division}\label{Janet}

In this subsection we will focus on the combinatorics of equations, enabling us to control the infinite set of differential variables appearing as partial derivatives of differential indeterminates.
For this purpose we use \textsc{Janet} division~\cite{GB1} which defines these combinatorics and also automatizes construction of integrability conditions. 
An overview of modern development can be found in \cite{GerI,Seiler} and the original ideas by {\sc Janet} are formulated in \cite{Janet}.
This is achieved by partitioning the set of differential variables into finitely many ``cones'' and ``free'' variables.
For creating this partition we present an algorithm for inserting new equations into an existing set of equations and adjusting the cone decomposition.
Apart from this insertion algorithm the only other adaptation of the algebraic \textsf{Decompose} algorithm \crossref{(}{algo_decompose}{)} will be the reduction algorithm presented in \textsection\ref{differential reduction}.

We fix a (differential) \textbf{ranking} $<$, which is defined as a total order on the differential variables such that
$u^{(k)}<\partial_j u^{(k)}$ and $u^{(k)}<u^{(l)}$ implies $\partial_j u^{(k)}<\partial_j u^{(l)}$ for all $u^{(k)},u^{(l)}\in U\mbox{, }\partial_j\in\Delta$.
For any finite set of differential variables, a differential ranking is a ranking as defined for the algebraic case in \crossref{\textsection}{algebraic_definition_notation}{}.
This allows us to define the largest differential variable $\ld(p)$ appearing in a differential polynomial $p\in F\{U\}$ as \textbf{leader}, which is set to $1$ for $p\in F$.
Furthermore, define $\rank(p)$ and $\ini(p)$ as the degree in the leader and the coefficient of $\ld(p)^{\rank(p)}$, respectively.
Again we will assume $1<u_{\mathbf{i}}^{(j)}$ for all $j\in\{1,\ldots,m\}$ and $\mathbf{i}\in\Z^n_{\ge0}$.

A set $W$ of differential variables is \textbf{closed} under the action of $\Delta'\subseteq\Delta$ if $\partial_i w\in W\ \forall\partial_i\in\Delta',w\in W$.
The smallest such closed set containing a differential variable $w$ denoted by $\langle w\rangle_{\Delta'}$ is called a \textbf{cone} and the elements of $\Delta'$ we shall call ({\sc Janet}) \textbf{admissible derivations}\footnote{ In~\cite{Ger3} and \cite[Chap.~7]{Seiler} the admissible derivations are called (\textsc{Janet}) multiplicative variables.}.
The $\Delta^{\prime}$-closed set generated by a set $W$ of differential variables is defined to be
\[
  \langle W\rangle_{\Delta^{\prime}}:=\bigcap_{\stackrel{W_i\supseteq W}{W_i\, \Delta^{\prime}\mbox{\scriptsize -closed}}} W_i\ \ \subseteq \ \ \langle U\rangle_{\Delta}\enspace.
\]

For a finite set $W=\{w_1,\ldots,w_r\}$, the {\sc Janet} \textbf{division} algorithmically assigns admissible derivations to the elements of $W$ such that the cones generated by the $w\in W$ are disjoint.
The derivation $\partial_l\in\Delta$ is  assigned to the cone generated by $w=u^{(j)}_{\mathbf{i}}\in W$ as admissible derivation, if and only if
\[
  \mathbf{i}_l=\max\left\{\mathbf{i}'_l\mid u^{(j)}_{\mathbf{i}'}\in W,\mathbf{i}'_k=\mathbf{i}_k\mbox{ for all }1\le k<l\right\}
\]
holds.
We remark, that $j$ is fixed in this definition, i.e., when constructing cones we only take into account other differential variables belonging to the same differential indeterminate.
The admissible derivations assigned to $w$ are denoted by $\Delta_W(w)\subseteq\Delta$ and we call the cone $\langle w\rangle_{\Delta_W(w)}$ the {\sc Janet} \textbf{cone} of $w$ with respect to $W$.
This construction ensures disjointness of cones but not necessarily that the union of cones equals $\langle W\rangle_\Delta$.
For the {\sc Janet} \textbf{completion} a finite set $\widetilde W\supset W$ is successively created by adding any 
$
  \tilde w=\partial_i w_j\not\in\biguplus_{w\in\widetilde W}\langle w\rangle_{\Delta_{\widetilde W}(w)}  
$
to $\widetilde W$, where $w_j\in\widetilde W$ and $\partial_i\in\Delta\setminus\Delta_{\widetilde W}(w_j)$. This leads to the disjoint {\sc Janet} \textbf{decomposition}
\[
  \langle W\rangle_\Delta=\biguplus_{w\in\widetilde W}\langle w\rangle_{\Delta_{\widetilde W}(w)}
\]
that separates a $\Delta$-closed set $W$ into finitely many cones $\langle w\rangle_{\Delta_{\widetilde W}(w)}$ after finitely many steps.
For details see \cite[Def.\ 3.4]{GerI} and \cite[Corr.\ 4.11]{GB1}.

With the {\sc Janet} decomposition being defined for sets of differential variables, we will assign admissible derivations to differential polynomials according to their leaders.
In particular, we extend the definitions of $\Delta_W(w)$ for finite $W\subset F\{U\}$ and $w\in W$ by defining $\Delta_W(w):=\Delta_{\ld(W)}(\ld(w))$.
  
A differential polynomial $q\in F\{U\}$ is called \textbf{reducible} with respect to  $p\in F\{U\}$, if there exists $\mathbf{i}\in\Z_{\ge 0}^n$ such that $\partial_1^{\mathbf{i}_1}\cdot\ldots\cdot\partial_n^{\mathbf{i}_n}\ld(p)=\ld(\partial_1^{\mathbf{i}_1}\cdot\ldots\cdot\partial_n^{\mathbf{i}_n}p)=\ld(q)$
and $\rank(\partial_1^{\mathbf{i}_1}\cdot\ldots\cdot\partial_n^{\mathbf{i}_n}p)\le\rank(q)$.
We call a derivative of an equation by an admissible derivation an \textbf{admissible prolongation}.
When restricting ourselves to admissible prolongation, we get the following concept:
For a finite set $T\subset F\{U\}$ we call a differential polynomial $q\in F\{U\}$ {\sc Janet} \textbf{reducible} with respect to  $p\in T$, if there exists $\mathbf{i}\in\Z_{\ge 0}^n$ 
such that $\partial_1^{\mathbf{i}_1}\cdot\ldots\cdot\partial_n^{\mathbf{i}_n}\ld(p)=\ld(q)$ with all proper derivatives being admissible
and $\rank(\partial_1^{\mathbf{i}_1}\cdot\ldots\cdot\partial_n^{\mathbf{i}_n}p)\le\rank(q)$ holds.
We shall also say that $q$ is {\sc Janet} \textbf{reducible} modulo $T$ if there is a $p\in T$, such that $q$ is {\sc Janet} reducible with respect to $p\in T$.

A set of differential variables $T\subset\langle U\rangle_{\Delta}$ is called \textbf{minimal}, if for any set $S\subset\langle U\rangle_{\Delta}$ with $\biguplus_{t\in T}\langle t\rangle_{\Delta_T(t)}=\biguplus_{s\in S}\langle s\rangle_{\Delta_S(s)}$ the condition $T\subseteq S$ holds (cf.\ \cite[Def.~4.2]{GB2}). We also call a set of differential polynomials minimal, if the corresponding set of leaders is minimal.

In addition to the non-zero initials and square-freeness of polynomials in the candidate set $S_T$ for a simple system (as defined in \crossref{\textsection}{algebraic_algorithms}{}), the equations in $(S_T)^=$ are required to have admissible derivations assigned to them.
When an equation $p$ is not reducible modulo $(S_T)^=$ it is added to $(S_T)^=$ and all polynomials in $S_T$ with a leader being a derivative of $\ld(p)$ are removed from $S_T$, ensuring minimality.
Furthermore, all non-admissible prolongations are created to be processed.
This is formulated in the following algorithm:

\begin{myalgorithm}[{\sf InsertEquation}]\ \\
\label{algo_differential_insertequation}
 \textit{Input:} A system $S'$ and a polynomial $p_=\in F\{U\}$ not reducible modulo $(S'_T)^=$.\\
 \textit{Output:} A system $S$, where $(S_T)^=\subseteq (S_T')^=\cup\{p_=\}$ is maximal satisfying
\[
\{\ld(q)\mid q\in(S_T)\setminus\{p\}\}\cap\langle\ld(p)\rangle_{\Delta}=\emptyset,
\]
\[
S_Q=S_Q'\cup(S'_T\setminus S_T)\cup\{(\partial_i q)_=\mid q\in (S_T)^=,\partial_i\not\in\Delta_{((S_T)^=)}(q)\}\enspace.
\]
 \textit{Algorithm:}
  \begin{algorithmic}[1]
    \STATE $S\gets S'$
    \STATE $S_T \gets S_T\cup\{p_=\}$
    \FOR{$q\in S_T\setminus\{p\}$}
      \IF{$\ld(q)\in\langle\ld(p)\rangle_{\Delta}$}
        \STATE $S_Q\gets S_Q\cup\{q\}$
        \STATE $S_T\gets S_T\setminus\{q\}$
     \ENDIF 
    \ENDFOR
    \STATE Reassign admissible derivations to $(S_T)^=$
    \STATE $S_Q\gets S_Q\cup \{(\partial_i q)_=\mid q\in (S_T)^=,\partial_i\notin\Delta_{((S_T)^=)}(q)\}$
    \RETURN $S$
  \end{algorithmic}
\end{myalgorithm}

\ifthenelse{\equal{\articletype}{split} \or \equal{\splitarticle}{jsc}}{The proof of correctness and termination is clear.}
We remark that a non-admissible prolongations might be added to $S_Q$ again each step, even though it has been added before.
This can be prevented by simply storing all previously generated non-admissible prolongations.

\subsection{Differential Simple Systems}\label{differential reduction}

This section goes on reducing the differential decomposition algorithm to the algebraic one.
We start by introducing partial solutions in order to algebraically evaluate differential polynomials at them yielding univariate differential polynomials.
Then we present a differential reduction algorithm, as the second distinction from the algebraic decomposition algorithm.
At last this section defines differential simple systems.

For a differential variable $x\in\langle U\rangle_{\Delta}$ and a power series $e\in E$ define the $F$-algebra homomorphism
\[
  \phi_{<x,e}: F\{U\}\to\overline{F}[\ v\mid v\in \langle U\rangle_{\Delta}, v\ge x\ ]:
    \begin{cases}
      u^{(j)}_{\mathbf{i}}\mapsto\alpha(e)(u^{(j)}), & \mbox{ for }u^{(j)}_{\mathbf{i}}<x\\
      u^{(j)}_{\mathbf{i}}\mapsto u^{(j)}_{\mathbf{i}}, & \mbox{ for }u^{(j)}_{\mathbf{i}}\ge x
    \end{cases}
\]
evaluating all differential variables  of a differential polynomial at $e$ which are $<$-smaller than $x$.

For differential reduction the {\sc Janet} partition of differential variables provides the mechanism to get a unique reductor in a fast way (for an algorithm see \cite{GYB}) by restricting to admissible prolongations.
After finding a reductor we apply a pseudo remainder algorithm (see\ \crossref{Eq.\ (}{prempquo}{)}).

We need to ensure that initials (and initials of the derivatives) of equations are non-zero.
Let $p\in F\{U\}$ with $x=\ld(p)$ and define the \textbf{separant} $\sep(p):=\frac{\partial p}{\partial x}$.
One easily checks (cf.\ \cite[\textsection I.8, lemma 5]{Kol} or \cite[\textsection3.1]{Hubert2}) that the initial of any derivative of $p$ is $\sep(p)$ and the separant of any square-free equation $p$ is nonzero on $\sol(p)$.
So by making the equations square-free, it is ensured that pseudo reductions are not only possible modulo $p$, but also modulo its derivatives.
This allows us to formulate the differential reduction algorithm:

\begin{myalgorithm}[{\sf Reduce}]\label{algo_differential_reduce}\ \\
 \textit{Input:} A differential system $S$ and a polynomial $p \in F\{U\}$.\\
 \textit{Output:} 
 A polynomial $q$ that is not \textsc{Janet} reducible modulo $S_T$ 
 with $\phi_e(p)=0$ if and only if $\phi_e(q)=0$ for each $e\in\sol(S)$.\\
 \textit{Algorithm:} \begin{algorithmic}[1]
                      \STATE $x \gets \ld(p)$
                       \WHILE{exists $q_=\in(S_T)^=$ and $i_1,\ldots,i_n\in\Z_{\ge0}$ with $i_j=0$ for $\partial_j\not\in\Delta_{(S_T)^=}(q)$ such that $\partial_1^{i_1}\cdot\ldots\cdot\partial_n^{i_n}\ld(q)=\ld(p)$ and  $\rank(\partial_1^{i_1}\cdot\ldots\cdot\partial_n^{i_n}p) \geq \rank(q)$ hold}
                       \STATE $p \gets \prem(p, \partial_1^{i_1}\cdot\ldots\cdot\partial_n^{i_n}q, x)$
                       \STATE $x \gets \ld(p)$
                      \ENDWHILE
                      \IF{$\textrm{\sf Reduce}(S, \ini(p)) = 0$}
                       \RETURN $\textrm{\sf Reduce}(S, p-\ini(p)x^{\rank(p)})$
                      \ELSE
                       \RETURN $p$
                      \ENDIF
                     \end{algorithmic}
\end{myalgorithm}

A polynomial $p\in F\{U\}$ is called \textbf{reduced\footnote{
There is a fine difference between not being reducible and being reduced.
In the case of not being reducible the initial of a polynomial can still reduce to zero and iteratively the entire polynomial.
}
modulo $S_T$} if $\mathsf{Reduce}(S, p)=p$.
A polynomial $p\in F\{U\}$ \textbf{reduces to $q$ modulo $S_T$} if $\mathsf{Reduce}(S, p)=q$.

\ifthenelse{\equal{\articletype}{split} \or \equal{\splitarticle}{jsc}}{%
The correctness is clear and termination is provided by \textsc{Dickson}'s Lemma (cf.\ \cite[Chap.~2, Thm.~5]{CLO} or \cite[\textsection 0.17, Lemma 15]{Kol}).
This lemma states that the ranking $<$ is well-founded on the set of leaders, i.e., a strictly $<$-descending chain of leaders is finite.
In our algorithm we have a descending chain of leaders $x$ of $p$ both in the loop and in the recursive call.
}%

Usually in differential algebra, one distinguishes a (full) differential reduction as used here and a partial (differential) reduction.
Partial reduction only employs \emph{proper} derivations of equations for reduction (cf.\ \cite[\textsection I.9]{Kol} or \cite[\textsection 3.2]{Hubert2}).
This is useful for separation of differential and algebraic parts of the algorithm and for the use of {\sc Rosenfeld}'s Lemma (cf.\ \cite{Rosenfeld}).

\begin{mydefinition}[Differential Simple Systems]
A differential system $S$ is (\textsc{Janet}) \textbf{involutive}, if 
all non-admissible prolongations in $(S_T)^=$ reduce to zero by $(S_T)^=$.\\
A system $S$ is called \textbf{differentially simple} or \textbf{simple}, if $S$ is
   \renewcommand{\theenumi}{\alph{enumi}}
   \renewcommand{\labelenumi}{\theenumi)}
\begin{enumerate}
  \item algebraically simple in the finite set of differential variables appearing in it,
  \item involutive,
  \item $S^=$ is minimal,
  \item no inequation is reducible modulo $S^=$.
\end{enumerate}
   \renewcommand{\theenumi}{\arabic{enumi}}
   \renewcommand{\labelenumi}{\theenumi)}
A disjoint decomposition of a system into differentially simple subsystems is called \textbf{(differential) }{\sc Thomas}\textbf{ decomposition}.
\end{mydefinition}

\subsection{The Differential Decomposition Algorithm}\label{differential algorithm}

The differential \textsc{Thomas} decomposition algorithm is a modification of the algebraic \textsc{Thomas} decomposition algorithm.
We have already introduced the new algorithms \textsf{InsertEquation} (\ref{algo_differential_insertequation}) for adding new equations into the systems and \textsf{Reduce} (\ref{algo_differential_reduce}) for reduction, that can replace their counterparts in the algebraic algorithm.

\begin{myalgorithm}[{\sf DifferentialDecompose}]\label{algo_differential_decompose}\ \\
 \textit{Input:} A differential system $S^\prime$ with $({S^\prime})_T=\emptyset$. \\
 \textit{Output:} A differential {\sc Thomas} decomposition of $S^\prime$. \\
 \textit{Algorithm:} The algorithm is obtained by replacing the two subalgorithms {\sf\nobreak InsertEquation} and {\sf\nobreak Reduce} in \crossref{(}{algo_decompose}{)} with their differential counterparts (\ref{algo_differential_insertequation}) and (\ref{algo_differential_reduce}), respectively.
 
\end{myalgorithm}

\ifthenelse{\equal{\articletype}{split} \or \equal{\splitarticle}{jsc}}{
\begin{proof}[Correctness]
We add three further loop invariants to be proven:
  
\begin{enumerate}
  \item $(S_T)^=$ is minimal.
  \item No inequation in $(S_T)^{\not=}$ is \textsc{Janet} reducible modulo $S_T$.
  \item For each system $S\in P\cup \mathit{Result}$ any non-admissible prolongation $p$ of $(S_T)^=$ reduces to zero using conventional differential reductions modulo $(S_Q)^=$ and using \textsc{Janet} reductions modulo $(S_T)^=$.
\end{enumerate}

The first loop invariant is a purely combinatorial matter, which is proven in \cite{Ger2} for an algorithm using exactly the same combinatorial approach.

Proving the second loop invariant is equally simple.
On the one hand, a newly added inequation $q$ in $S_T$ is not \textsc{Janet} reducible by $(S_T)^=$, since algorithm \textsf{Reduce} (\ref{algo_differential_reduce}) is applied to it before.
On the other hand, algorithm {\sf InsertEquation} (\ref{algo_differential_insertequation}) removes all inequations from $S_T$ being divisible by a newly added equation and places them into $S_Q$.

The non-trivial third loop invariant clearly holds at the beginning of the algorithm, because $S_T$ is empty.

We claim that reduction of an equation $q_=\in S_Q$ by $(S_T)^=$ in \crossref{line }{reduce_q}{ of Algorithm \crossref{(}{algo_decompose}{)}} does not effect the loop invariant.
We prove this claim for single reductions, which generalizes by an easy induction.
Let $q':=\prem(q,p,x)=m\cdot q-\pquo(q,p,x)\cdot p$ be a pseudo remainder equation
(see \crossref{(}{prempquo}{) on page \pageref{prempquo}})
reducing $q$ to $q'$.
Then a pseudo remainder equation $\prem(r,q,x)=m'\cdot r-\pquo(r,q,x)\cdot q$ describing a reduction by $q$ might simply be rewritten as an iterated pseudo remainder equation
\[
  \underbrace{m\cdot \prem(r,q,x)}_{\prem(\prem(r,p,x),q',x)}=m\cdot 
  m'\cdot r-\pquo(r,q,x)\cdot q'-\pquo(r,q,x)\cdot \pquo(q,p,x)\cdot p
\]
and using the \textsc{Leibniz} rule the same holds for reduction modulo partial derivatives of $q$.
So a reduction to zero modulo $q$ and its derivatives can be achieved by iterated pseudo reduction modulo $q'$, $p$, and their derivatives.

This holds especially for an equation $q_=\in S_Q$ reducing to $0$ modulo $(S_T)^=$ in \crossref{line }{reduce_q}{}, which can be removed from $S_Q$ without violating the loop invariant.

When the square-free part of $p_=$ of $q_=$ is inserted into $S_T$ in \crossref{line }{q_2}{}, all non-admissible prolongations of $p_=$ are added to $S_Q$ as equations and thus reduce to $0$ modulo $(S_Q)^=$.
In addition, moving equations from $S_T$ back into $S_Q$ when inserting a new one does not change the loop invariant either, because afterwards reductions are allowed using these equations and \emph{all} their prolongations in contrast to admissible prolongationss \emph{only}, as allowed before.
Furthermore, every non-admissible prolongation reduced to zero using $q_=\in (S_Q)^=$ now \textsc{Janet} reduces to zero modulo $p_=$ in $S_T$ or conventionally reduces to zero modulo its non-admissible prolongations inserted into $S_Q$.
This holds for two reasons:
First, the non-admissible prolongations are exactly the derivatives of $p_=$, which are not allowed to be used when applying \textsc{Janet} reduction.
Second, write $m\cdot q=p\cdot q_1$ with $\ld(m)<x$ and $\pma(m)\not=0\ \forall\mathbf{a}\in\sol(S_{<\ld(q)})$.
Then $p$ as a divisor of $q$ algebraically reduces $q$ to zero and the derivative $\partial q$ of $q$ is implied by $p_=$ and $(\partial p)_=$, since $\partial(m\cdot q)=(\partial p)\cdot q_1+p\cdot(\partial q_1)$ for any $\partial\in\Delta$.
Inductively, a repeated derivative of $q$ is implied by the repeated derivative of $p$.

When computing the gcd of two equations in \crossref{line }{resplitgcd}{}, the gcd of $q$ and $(S_T)_x$ will be inserted into $S_T$ and reduces everything to zero that both $q$ and $(S_T)_x$ did.
As above, the non-admissible prolongations are covered by inserting them into $S_Q$ and the admissible prolongations are implied.
In \crossref{line}{insertres}{} where the computation of the gcd is prepared, the equation $q_=$ will simply be reinserted into $S_Q$, being already reduced.

Dividing an equation $(S_T)_x$ by an inequation $q_{\not=}$ in \crossref{lines }{dividebyinequation}{} and \crossref{}{q_3}{} also influences $(S_T)^=$.
The new equation $p$, being a divisor of $(S_T)_x$, reduces everything to zero that $(S_T)_x$ and its non-admissible prolongations did by the above arguments.
This proves the third loop invariant.

When the algorithm terminates, $S_Q$ is empty and thus all non-admissible prolongations from $(S_T)^=$ \textsc{Janet} reduce to zero by $(S_T)^=$ by the third loop invariant.
Thus the system is involutive.

One easily checks that the correctness proof of the algebraic decomposition algorithm {\sf Decompose} applies to the differential case, since we use power series as a solution set.
Thus we see there that each system in the output is algebraically simple.

Furthermore, the first loop invariant directly ensures minimality and the second loop invariant ensures that no inequation is reducible by an equation, since for an involutive set reducibility modulo this set is equivalent to \textsc{Janet} reducibility.
Combining these three facts, the result of the differential decomposition is also a disjoint decomposition algorithm into differentially simple systems.\qed
\end{proof}
}{}

\ifthenelse{\equal{\articletype}{split} \or \equal{\splitarticle}{jsc}}{
Our main tool for proving the termination of the algorithm is again using six orders on systems.
To show that the orders, used in addition to the algebraic ones introduced in \crossref{(}{def_orders}{)}, are well-founded we use {\sc Dickson}'s lemma again.

\begin{defrem} We define the orders $\prec_{1a}$, $\prec_{1b}$, $\prec_{1c}$, $\prec_{2}$, $\prec_{3}$, and $\prec_{4}$ on systems and show their well-foundedness.
Using \crossref{Remark (}{composite_well_founded}{)}, $\prec:=[{\prec_{1a},} {\prec_{1b},} {\prec_{1c},} {\prec_{2},} {\prec_{3},} {\prec_{4}}]$ is again well-founded as the composite order.
\begin{itemize}
  \item[$\prec_{1a}$:] For $V\subseteq\langle U\rangle_{\Delta}$ there is a unique minimal set $\nu(V)\subseteq V$ with $V\subseteq \langle \nu(V)\rangle_{\Delta}$ (cf.\ \cite[Chap.~2, \textsection4, exercise 7 and 8]{CLO}), called \textbf{canonical differential generators} of $V$.
  For a system $S$, define $\nu(S)$ as $\nu(\ld((S_T)^=))$.\\  
  For systems $S, S'$ we define $S\prec_{1a}S'$ if and only if 
 $\min_<(\nu(S)\setminus\nu(S'))<\min_<(\nu(S')\setminus\nu(S))$. An empty set is assumed to have $x_\infty$ as minimum, which is $<$-larger than all other leaders.\\
  By {\sc Dickson}'s lemma, $\prec_{1a}$ is well-founded.

 \item[$\prec_{1b}$:] 
 For systems $S, S'$ define $S\prec_{1b}S'$ if and only if  $\min_<\left(\ld(S)\setminus\ld(S')\right)<\min_<\left(\ld(S')\setminus\ld(S)\right)$ and $S\osim_{1a} S'$ .
 An empty set is assumed to have $x_\infty$ as minimum, which is $<$-larger than all other leaders.\\
  Minimality of $(S_T)^=$ at each step of the algorithm and the constructivity property of the {\sc Janet} division (cf.\ \cite[Prop.~4.13]{GB1}) imply well-foundedness of $\prec_{1b}$ (cf.\ \cite[Thm.~4.14]{GB1}).

 \item[$\prec_{1c}$:] 
 For systems $S$ and $S'$ with $S\osim_{1a}S'$ and $S\osim_{1b}S'$ both having leaders $x_1,\ldots,x_l$ in $(S_T)^=$ and $(S'_T)^=$ define $S\prec_{1c,x_k} S'$ if and only if $\rank((S_T)^=_{x_i})<\rank((S'_T)^=_{x_i})$. This ordering is clearly well-founded.
 For these systems define $S\prec_{1c}S'$ as $[\prec_{1c,x_1},\ldots,\prec_{1c,x_l}]$, which as a composite ordering is again well-founded.

 \item[$\prec_2$:] This is defined almost identical to the algebraic $\prec_2$. In this case the set of possible leaders is $\{1\}\cup\langle U\rangle_{\Delta}$.\\
 To show well-foundedness of the differential ordering $\prec_2$ we use that $<$ is well-founded on the set of leaders as implied by {\sc Dickson}'s lemma.
 This way, $<$ is extended to a well-founded ordering on $\operatorname{im}(\mu)=\{1,x_{\infty}\}\cup\langle U\rangle_{\Delta}$ by $1<x_{\infty}$ and $y<x_{\infty}$ for all $y\in\langle U\rangle_{\Delta}$.

 \item[$\prec_3$:] This is verbatim the same condition and proof of well-foundedness as in the algebraic case.
 But we remark that we do a {\sc Noether}ian induction (cf.\ \cite[III.6.5, Prop.~7]{Bou68}) instead of an ordinary induction.

 \item[$\prec_4$:] This is identical to the algebraic case.
\end{itemize}
\end{defrem}

\begin{proof}[Termination]

Showing that systems get $\prec$-smaller is a trivial task.
This is clear in the case of $\prec_2$, $\prec_3$, and $\prec_4$, since these properties can be used equivalently for the algebraic and differential algorithm.
In the algebraic case a system $\prec_1$-decreases if and only if either an equation is added to $S_T$ or the degree of an existing equation in $S_T$ is decreased.
Now, in the differential case, inserting a completely new equation decreases either $\prec_{1a}$ or $\prec_{1b}$.
If the degree of a polynomial gets lower, the system $\prec_{1c}$-decreases.

Thus, like in the algebraic termination proof, we have an asymmetric, transitive, and well-founded order on systems, which decreases in each step of the loop. This proves termination.\qed
\end{proof}
}{}

We give an example taken from \cite[pp.\ 597-600]{BC}:

\begin{myexample}[Cole-Hopf Transformation]
For $F:=\R(x,t)$, $\Delta=\{\frac{\partial}{\partial x},\frac{\partial}{\partial t}\}$, and $U=\{\eta,\zeta\}$ consider the heat equation $h=\eta_t+\eta_{xx}\in F\{U\}^=$ and {\sc Burger}'s equation $b=\zeta_t+\zeta_{xx}+2\zeta_x\cdot\zeta\in F\{U\}^=$.
To improve readability, leaders of polynomials are underlined below.

First we claim that any power series solution for the heat equation with a non-zero constant term can be transformed to a solution of \textsc{Burger}'s equation by means of the {\sc Cole-Hopf} transformation $\lambda: \eta\mapsto\frac{\eta_x}{\eta}$.
The differential {\sc Thomas} decomposition for an orderly ranking with $\zeta_x>\eta_t$ of
\[
  \{h_=,\underbrace{(\eta\cdot\zeta-\eta_x)_=}_{\Leftrightarrow \zeta=\lambda(\eta)},\eta_{\not=}\}
\]
consists of the single system
\[
  S=\{(\underline{\eta_x}-\eta\cdot\zeta)_=,(\eta\cdot\underline{\zeta_x}+\eta_t+\eta\cdot\zeta^2)_=,\underline{\eta}_{\not=}\}
\]
and one checks that $\mathsf{Reduce}(S,b)=0$ holds. This implies that any non-zero solution of the heat equation is mapped by the {\sc Cole-Hopf} transformation to a solution of {\sc Burger}'s equation.

In addition we claim that $\lambda$ is surjective.
For the proof we choose an elimination ranking (cf. \cite[\textsection 8.1]{Hubert2} or \cite{BoulierElimination}) with $\eta\gg\zeta$, i.e., $\eta_\mathbf{i}>\zeta_\mathbf{j}$ for all $\mathbf{i},\mathbf{j}\in \left(\mathbb{Z}_{\ge 0}\right)^2$.
We compute the differential {\sc Thomas} decomposition of $\{h_=,b_=,(\eta\cdot\zeta-\eta_x)_=,\eta_{\not=}\}$ which again consists of a single system
\[
  S=\{(\underline{\eta_x}-\eta\cdot\zeta)_=,(\eta\cdot\zeta_x+\underline{\eta_t}+\eta\cdot\zeta^2)_=,\underline{b}_=,\underline{\zeta}_{\not=}\}\enspace.
\]
The properties of a simple system ensure that for any solution of lower ranking equations there exists a solution of the other equations (cf.\ \crossref{(}{exist_sol}{)}).
The elimination ordering guarantees that the only constraint for $\zeta$ is {\sc Burger}'s equation $b_=$ and thus for any solution $f\in\sol(b_=)$ there exists a solution $(g,f)\in\sol(S)$.
Furthermore, since $h_=$ was added to the input system, $g\in\sol(h_=)$ holds and finally the equation $(\eta\cdot\zeta-\eta_x)_=$ implies $\lambda(g)=f$.
\end{myexample}

\begin{myremark}
Elements of the differential field are not subjected to splittings, unless they are modelled as differential indeterminates.
For example to model a differential field $F=\C(x)$ with $\Delta=\{\frac{\partial}{\partial x}\}$, we add an extra differential indeterminate $X$ to $U$ and replace $x$ by $X$ in all equations and inequations.
We subject $X$ to the relation $\frac{\partial}{\partial x}X=1$ for $X$ being ``generic'' or  $(\frac{\partial}{\partial x}X-1)\cdot\frac{\partial}{\partial x}X=0$, if we allow $X$ to degenerate to a point.
This will be subject of further study.
\end{myremark}

%% file: implementation.tex
\section{Implementation}\label{implementation}

\subsection{Algorithmic Optimizations}

\ifthenelse{\equal{\articletype}{combined} \or \equal{\splitarticle}{algebraic}}{%
In the \textsf{Decompose} algorithm, pseudo remainder sequences for the same pairs of polynomials are usually needed several times.
As these calculations are expensive in general, for \emph{avoiding repeated calculations}, it is important that the results are kept in memory and will be reused when the same sequence is requested again.

If a polynomial admits \emph{factorization}, we can use the it to save computation time.
More precisely, a disjoint decomposition of the system $S \uplus \{ (p\cdot q)_= \}$ is given by $(S \cup \{ p_= \}, S \cup \{ p_{\neq}, q_= \})$ and the system $S \uplus \{ (p\cdot q)_{\neq} \}$ is equivalent to $S \cup \{ p_{\neq}, q_{\neq} \}$.
Let $Y_i:=\left\{x_j\mid x_j<x_i,(S_T)^=_{x_j}\not=\emptyset\right\}$ and $Z_i:=\left\{x_j\mid x_j<x_i, (S_T)^=_{x_j}=\emptyset\right\}$.
If $(S_T)^=_{x_i}$ is irreducible over the field $F_i:=F(Z_i)[Y_i]/\langle (S_T)^=_{<x_i}\rangle_{F(Z_i)[Y_i]}$ for all $i\in\{1,\ldots,n\}$, where $\langle (S_T)^=_{<x_i}\rangle_{F(Z_i)[Y_i]}$ is the ideal generated by $(S_T)^=_{<x_i}$ in the polynomial ring $F(Z_i)[Y_i]$, factorization of polynomials can be performed over $F_n$ instead of $F$.

Coefficient growth is a common problem in elimination.
If possible, polynomials should be represented as compact as possible.
Once it is known that the initial cannot vanish, the \emph{content} (in the univariate sense) cannot vanish either.
Thus, every time an initial has been added as an inequation to the system, one can divide the polynomial by its content.

If the ground field $F$ is represented as a field of fractions of a domain $D$ (like the rationals or a rational function field over the rationals), it also makes sense to remove the multivariate content, which is an element of $F$.

When reducing, in addition to reduction modulo the polynomial of the same leader, reducing the coefficients modulo the polynomials of lower leader can be considered.
In some cases this leads to a reduction of sizes of coefficients, in other cases sizes increase.
The latter is partly due to whole polynomials being multiplied with initials of the reductors.
Finding a good heuristic for this \emph{coefficient reduction} is crucial for efficiency.

In the algebraic algorithm, polynomials don't necessarily have to be square-free when they are inserted into the candidate simple system.
Efficiency is sometimes improved greatly by postponing the calculation of the square-free split as long as possible.
}{}%

\ifthenelse{\equal{\splitarticle}{differential}}{
Most algebraic optimizations as mentioned in \crossref{}{BGLROptimizations}{} are also valid in the differential case, with the exception of not being able to forgo making polynomials squarefree from the beginning (cf.\ \ref{differential reduction}). 
}{}%
\ifthenelse{\equal{\articletype}{combined}}{In the differential case application }{%
\ifthenelse{\equal{\splitarticle}{differential}}{Furthermore, application }{}}%
\ifthenelse{\equal{\articletype}{combined} \or \equal{\splitarticle}{differential}}{%
of \emph{criteria} to avoid useless reduction of non-admissible prolongations can decrease computation time.
The combinatorial approach used in this paper already avoids many reductions of so-called $\Delta$-polynomials, as used in other approaches (see \cite{GY}).
Nonetheless, using the involutive criteria 2-4 (cf. \cite{GB1,GerI,AH} and \cite[\textsection 4, Prop.\ 5]{BLOP}) which together are equivalent to the chain criterion, is valid and helpful.%
}{}%

\ifthenelse{\equal{\articletype}{combined} \or \equal{\splitarticle}{algebraic}}{%
Another possible improvement is parallelization, since the main loop in line \ref{main_loop} of \textsf{Decompose} (\ref{algo_decompose}) can naturally be used in parallel for different systems.%
}{}%

\subsection{Implementation in {\sc Maple}}

\ifthenelse{\equal{\articletype}{combined}}{Both algorithms have }{The algorithm has }%
been implemented in the \textsc{Maple} computer algebra system.
Packages can be downloaded from \cite{hp_thomasdecomp}, documentation and example worksheets are available there.

\ifthenelse{\equal{\articletype}{combined} \or \equal{\splitarticle}{algebraic}}{%
The main reason for choosing \textsc{Maple} for the implementation is the collection of solvers for polynomial equations, ODEs, and PDEs already present.
Furthermore, fast algorithms exist for polynomial factorization over finitely generated field extensions of $\mathbb{Q}$ and for gcd computation.
Computation of subresultants is not available in \textsc{Maple}, therefore an algorithm based on \cite{Ducos} is implemented for that purpose.%
}{}

\ifthenelse{\equal{\articletype}{combined} \or \equal{\splitarticle}{differential}}{%
Features for the differential package include arbitrary differential rankings, using functions implemented in \textsc{Maple} as differential field, computation of power series solutions, and a direct connection to the solvers of \textsc{Maple} for differential equations.

\input{implementation_example_diff}}{}

\subsection{Implementations of Similar Decomposition Algorithms}

\ifthenelse{\equal{\articletype}{combined} \or \equal{\splitarticle}{algebraic}}{%
The \textsf{RegularChains} package \cite{regularchains}, which is shipped with recent versions of \textsc{Maple}, implements a decomposition of a polynomial ideal into ideals represented by regular chains and a radical decomposition of an ideal into square-free regular chains.
The solution sets of this decomposition are in general not disjoint.
However, there is an extension called comprehensive triangular decomposition (cf.~\cite{ctd}) that provides disjointness on the parameters of a parametric system.
The systems of the parameters are not simple systems though.
The \textsf{RegularChains} package contains \textsf{FastArithmeticTools} as a subpackage implementing asymptotically fast polynomial arithmetic for the modular case.

The $\epsilon$\textsf{psilon} package (\cite{epsilon}) by Dongming Wang implements different kinds of triangular decompositions (including a decomposition into regular chains like above) in \textsc{Maple}.
It is the only software package besides our own that implements the \textsc{Thomas} decomposition.
It uses the simpler ``top-down'' approach that Thomas (cf.\ \cite{Tho1,Tho2}) suggested, i.e., polynomials of higher leader are considered first.
All polynomials of the same leader are combined into one common consequence.
New systems, which contain conditions on initials of polynomials and subresultants, are created by splitting subalgorithms similar to ours.
All these new conditions of lower leader are not taken into account for now and will be treated in a later step.
Contrary to our approach, one cannot reduce modulo an \emph{unfinished} system and hence inconsistency checks are less natural and more complicated.
It is conceivable that this strategy spends too much time on computations with inconsistent systems.
Therefore, $\epsilon$\textsf{psilon} implements highly sophisticated heuristics for early detection of inconsistent systems.
It achieves similar performance to our implementation.%
}{%
The $\epsilon$\textsf{psilon} package (\cite{epsilon}) by Dongming Wang mentioned \cite[\textsection3]{BGLR1} also computes a \textsc{Thomas} decomposition for the ordinary differential case, but does not treat partial differential equations.%
}

\ifthenelse{\equal{\articletype}{combined} \or \equal{\splitarticle}{differential}}{%
The \textsc{Maple} package \textsf{diffalg} \cite{diffalg} deals with ordinary and partial differential equations as described in \cite{BLOP}.
Its functionalities are used by symbolic differential equations solvers in \textsc{Maple}.
For an input of equations and inequations it computes a radical decomposition of the differential ideal generated by the equations and saturated by the inequations.
I.e., a description of the vanishing ideal of the \textsc{Kolchin} closure (cf.\ \cite[\textsection IV.1]{Kol}) of the solutions is computed.
The output are differential characteristic sets as introduced by \textsc{Ritt} \cite[\textsection I.5]{Ritt}.
Computation of differential consequences is driven by reduction of $\Delta$-polynomials, which are the analogon of $s$-polynomials in differential algebra.
We found the system being optimized and well-suited for computations with ordinary differential equations.

Similar algorithms as in diffalg are used in the \textsf{BLAD}-libraries \cite{blad}.
It is designed as a stand-alone \textsf{C}-library with an emphasis on usability for non-mathematicians and extensive documentation.
As it is written in \textsf{C}, \textsf{BLAD} is expected to outperform \textsf{diffalg} for relevant examples.%
}{}

For future publications, we plan to compare the \textsc{Thomas} decomposition and our implementation with other decompositions and implementations.
We also plan to further examine applications that benefit from the properties of simple systems.

%% file: implementation_example_diff.tex
\begin{myexample}

\begin{maplegroup}
Start by loading the current version of our package:

\end{maplegroup}
\begin{maplegroup}
\emptyline
\end{maplegroup}
\begin{maplegroup}
\begin{mapleinput}
\mapleinline{active}{1d}{with(DifferentialThomas):
ComputeRanking([t],[x2,x1,y,u],"EliminateFunction");}{%
}
\end{mapleinput}

\end{maplegroup}
\begin{maplegroup}
\emptyline
\end{maplegroup}
\begin{maplegroup}
This creates the differential polynomial ring
$\mathbb{Q}\{x^{(2)},x^{(1)},y,u\}$ for
$\Delta=\{\frac{\partial}{\partial t}\}$. Here $u$ indicates the
input, $x^{(1)}$ and $x^{(2)}$ the state, and $y$ the output of the
system. The chosen ranking ``$<$'' is the elimination ranking with
$x^{(2)}\gg x^{(1)}\gg y\gg u$, i.e.,
$x^{(2)}_\mathbf{i}>x^{(1)}_\mathbf{j}>y_\mathbf{k}>u_\mathbf{l}$ for
all $\mathbf{i},\mathbf{j},\mathbf{k},\mathbf{l}\in \mathbb{Z}_{\ge
0}$.

\end{maplegroup}
\begin{maplegroup}
\emptyline
\end{maplegroup}
\begin{maplegroup}
\begin{mapleinput}
\mapleinline{active}{1d}{L:=[x1[1]-u[0]*x2[0],x2[1]-x1[0]-u[0]*x2[0],y[0]-x1[0]];}{%
}
\end{mapleinput}

\mapleresult
\begin{maplelatex}
\mapleinline{inert}{2d}{L := [x1[1]-u[0]*x2[0], x2[1]-x1[0]-u[0]*x2[0], y[0]-x1[0]];}{%
\[
L := [{\mathit{x1}_{1}} - {u_{0}}\,{\mathit{x2}_{0}}, \,{\mathit{
x2}_{1}} - {\mathit{x1}_{0}} - {u_{0}}\,{\mathit{x2}_{0}}, \,{y_{
0}} - {\mathit{x1}_{0}}]
\]
}
\end{maplelatex}

\end{maplegroup}
\begin{maplegroup}
\emptyline
\end{maplegroup}
\begin{maplegroup}
We follow \cite[Ex.~1]{diop} and want to compute the external
trajectories of a differential ideal generated by $L$, i.e. intersect
this differential ideal with $\mathbb{Q}\{y,u\}$. 

\end{maplegroup}
\begin{maplegroup}
\emptyline
\end{maplegroup}
\begin{maplegroup}
\begin{mapleinput}
\mapleinline{active}{1d}{res:=DifferentialThomasDecomposition(L,[]);}{%
}
\end{mapleinput}

\mapleresult
\begin{maplelatex}
\mapleinline{inert}{2d}{res := [DifferentialSystem, DifferentialSystem];}{%
\[
\mathit{res} := [\mathit{DifferentialSystem}, \,\mathit{
DifferentialSystem}]
\]
}
\end{maplelatex}

\end{maplegroup}
\begin{maplegroup}
\emptyline
\end{maplegroup}
\begin{maplegroup}
We show the equations and inequations of the differential systems not
involving $x^{(1)}$ and not involving $x^{(2)}$.
The chosen ranking guarantees that for each differential system of the
output, all constraints holding for lower ranking differential
indeterminates can be read off the equations and inequations only
involving these differential indeterminates, i.e., the systems shown
determine the external trajectories of the system:

\end{maplegroup}
\begin{maplegroup}
\emptyline
\end{maplegroup}
\begin{maplegroup}
\begin{mapleinput}
\mapleinline{active}{1d}{PrettyPrintDifferentialSystem(res[1]):
remove(a->has(a,x2) or has(a,x1),\%);}{%
}
\end{mapleinput}

\mapleresult
\begin{maplelatex}
\mapleinline{inert}{2d}{[-u(t)*diff(y(t),`$`(t,2))+diff(y(t),t)*u(t)^2+diff(y(t),t)*diff(u(t)
,t)+y(t)*u(t)^2 = 0, u(t) <> 0];}{%
\[
[ - \mathrm{u}(t)\,({\frac {d^{2}}{dt^{2}}}\,\mathrm{y}(t)) + (
{\frac {d}{dt}}\,\mathrm{y}(t))\,\mathrm{u}(t)^{2} + ({\frac {d}{
dt}}\,\mathrm{y}(t))\,({\frac {d}{dt}}\,\mathrm{u}(t)) + \mathrm{
y}(t)\,\mathrm{u}(t)^{2}=0, \,\mathrm{u}(t)\neq 0]
\]
}
\end{maplelatex}

\end{maplegroup}
\begin{maplegroup}
\begin{mapleinput}
\mapleinline{active}{1d}{PrettyPrintDifferentialSystem(res[2]):
remove(a->has(a,x2) or has(a,x1),\%);}{%
}
\end{mapleinput}

\mapleresult
\begin{maplelatex}
\mapleinline{inert}{2d}{[diff(y(t),t) = 0, u(t) = 0];}{%
\[
[{\frac {d}{dt}}\,\mathrm{y}(t)=0, \,\mathrm{u}(t)=0]
\]
}
\end{maplelatex}

\end{maplegroup}
\begin{maplegroup}
\emptyline
\end{maplegroup}
\begin{maplegroup}
These systems, having disjoint solution sets, are identical to the
ones found in \cite{diop}.

\end{maplegroup}

\end{myexample}

%% file: acknowledgements.tex
\section{Acknowledgements}\label{section_acknowledgements}

The contents of this paper profited very much from numerous useful comments and remarks from Wilhelm Plesken.
The authors thank him as well as Dongming Wang and Fran{\c{c}}ois Boulier for fruitful discussions.
Furthermore, our gratitude goes to the anonymous referees for valuable comments and for pointing out informative references.
The second author (V.P.G.) acknowledges the Deutsche Forschungsgemeinschaft for the financial support that made his stay in Aachen possible. The presented results were obtained during his visits.